\documentclass[11pt]{article}
\usepackage{amsthm,amsmath,amssymb,anysize}
\newtheorem{lemma}{Lemma}[section]
\newtheorem{theorem}[lemma]{Theorem}
\newtheorem{remark}[lemma]{Remark}
\newtheorem{proposition}[lemma]{Proposition}

\newtheorem{corollary}[lemma]{Corollary}
\marginsize{3cm}{3cm}{2.3cm}{2.2cm}%×óÓÒÉÏÏÂ
\numberwithin{equation}{section}

\begin{document}
\title{\textbf{Derivations for the even part of the
Hamiltonian superalgebra in positive characteristic\footnote{Supported by a NSF grant 10471096 of China and
``One Hundred Talents Program" from University of Science and Technology of China }}}
\author{Wende Liu$^{1,2,\,}$\footnote{Corresponding author. Email-addresses:
wendeliu@ustc.edu.cn (W. Liu), ycsu@ustc.edu.cn (Y. Su), zhyz@nenu.edu.cn (Y. Zhang)},
Yucai Su$^{1}$ and Yongzheng Zhang$^{2}$  \\
\\
{$^{1}$Department of Mathematics,} \\
{ University of Science and Technology of China, Hefei 230026,   China}\\
\\
{$^{2}$Department of Mathematics,} \\
{Harbin Normal University, Harbin 150080,   China}
}
\date{ }
\maketitle
\begin{quotation}
\small\noindent \textbf{Abstract}: \noindent In this paper we
consider the derivations for even part of the finite-dimensional
Hamiltonian superalgebra $H$ over a field of prime characteristic.
We first introduce an  ideal $\frak{N}$ of $H_{\overline{0}} $
and  show that   the derivation space from $H_{\overline{0}}$ into
$W_{\overline{0}}$ can be obtained by the derivation space from
$\frak{N}$ into $W_{\overline{0}},$ the even part of the
generalized Witt superalgebra $W$. For further application we also
give the generating set of the  ideal $\frak{N}$. Then
we describe three series of exceptional derivations from
$ {H}_{\overline{0}}$ into $W_{\overline{0}}.$ Finally, we determine all the derivations vanishing
on the non-positive $\mathbb{Z}$-graded part of $ H _{\overline{0}}$, the odd
$\mathbb{Z}$-homogeneous derivations, and negative
$\mathbb{Z}$-homogeneous derivations from  $ {H}_{\overline{0}} $ into
$W_{\overline{0}}.$
\\

\noindent  Mathematics Subject Classification 2000: 17B50, 17B40
\end{quotation}\vspace{0.5cm}

\noindent {\Large\textbf{0. \ \  Introduction}}
\newline

 \noindent
During last few decades  the theory of Lie superalgebras has
undergone a remarkable evolution both  in mathematics and in
physics (see \cite{Sch1}). For example,
  the classifications by V.G. Kac of
finite-dimensional simple Lie superalgebras and
infinite-dimensional simple linearly compact Lie superalgebras
over algebraically closed fields  of characteristic zero  have been
completed (see \cite{Kac2}, \cite{Kac4}). For modular Lie
superalgebras,
 as far as we
 know, \cite{KL} and  \cite{P} may be the earliest papers.

In this paper we consider derivations for the even parts of  modular
 Lie superalgebras of Cartan type $H.$
Our work is  originally  motivated by the work on modular Lie algebras of Cartan type
(see \cite{Ce,st,SF}).  Note that  the superderivation algebras have been
determined for the finite-dimensional
  modular Lie superalgebras of
 Cartan type
$W,$ $S,$ $H,$ and $K.$  (see \cite{MZ,WZ3,ZZ1}). The superderivation
algebra was determined for the finite-dimensional odd Hamiltonian superalgebra
$HO$ in \cite{LZ1}. We should mention that the
derivations of the even
parts have also been studied sufficiently for the Lie superalgebras of Cartan type $W,$  $S$
and $HO$ (see \cite{lhs,lz1});
in particular, the derivations from the even parts into the odd parts have been determined for
$W$ and $S$ (see \cite{lg}). However, the present work differs greatly from the ones mentioned above; in particular,
we find ``more" outer derivations for the even part of $H$ than $HO$, $W,$ or $S $ (see also Remark \ref{r3.11}).

This paper is organized as follows. In Section 1 we give the necessary notation and concepts. In Section 2,
we first introduce an ideal $\frak{N}$ of the even part of $H,$ which is crucial for our aim. Then we give the
generating set of the ideal $\frak{N}$ for future application. In Section 3 we mainly construct three series of outer
derivations from the even part of $H$ into the even part of $W,$ the generalized Witt superalgebra. In Section 4,
we determine all the derivations vanishing
on the non-positive $\mathbb{Z}$-graded part of $ H _{\overline{0}}$, the odd
$\mathbb{Z}$-homogeneous derivations, and negative
$\mathbb{Z}$-homogeneous derivations from  $ {H}_{\overline{0}} $ into
$W_{\overline{0}}.$

\section{Preliminaries}

   Let $\mathbb{Z}_{2}=\{\overline{0}, \overline{1}\}$ be the field
of two elements. For a   vector superspace
$V=V_{\overline{0}}\oplus V_{\overline{1}}, $ we denote by
$\mathrm{p}(a)=\theta$ the \textit{parity of a homogeneous
element} $a\in V_{\theta}, \theta\in \mathbb{Z}_{2}.$
 We assume throughout  that the notation   $\mathrm{p}(x)$ implies
that
  $x$ is a $\mathbb{Z}_2$-homogeneous element.

Let $\frak{g}$ be a  Lie algebra and $V$   a $\frak{g}$-module. A
linear mapping $D:\frak{g}\rightarrow V$ is called a
\textit{derivation } from $\frak{g}$  into $V$ if $D(xy)=x\cdot
D(y)-y\cdot D(x)$ for all $x,y\in \frak{g}.$ A derivation
$D:\frak{g}\rightarrow V$ is called \textit{inner} if there is
$v\in V $ such that $D(x)=x\cdot v$ for all $x\in \frak{g}.$
Following \cite[p. 13]{SF}, denote by $\mathrm{Der}(\frak{g},V)$
the \textit{derivation space} from $\frak{g}$  into $V.$  Then
$\mathrm{Der}(\frak{g},V)$ is a $\frak{g}$-submodule of
$\mathrm{Hom}_{\mathbb{F}}(\frak{g},V).$ Assume in addition that
$\frak{g}$ and $V$ are finite-dimensional and that
$\frak{g}=\oplus_{r\in \mathbb{Z}} \frak{g}_{_{[r]}}$ is
$\mathbb{Z}$-graded and $V=\oplus_{r\in\mathbb{ Z}} V_{[r]}$ is a
$\mathbb{Z}$-graded $\frak{g}$-module. Then
$\mathrm{Der}(\frak{g},V)=\oplus_{r\in \mathbb{Z}}
\mathrm{Der}_{[r]}(\frak{g},V)$ is a $\mathbb{Z}$-graded
$\frak{g}$-module  by setting
$$\mathrm{Der}_{[r]}(\frak{g},V):=\{D\in
\mathrm{Der}(\frak{g},V)\mid D(\frak{g}_{_{[i]}})\subset V _{[r+i]} \
\mbox{for all} \ i \in\mathbb{Z}\}.$$
In the case $V=\frak{g},$
the \textit{derivation algebra} $\mathrm{Der}(\frak{g})$ coincides
with $\mathrm{Der}(\frak{g},\frak{g}) $ and
$\mathrm{Der}(\frak{g}) =\oplus _{r\in \mathbb{Z}}
 \mathrm {Der}_{[r]}(\frak{g})$ is a $\mathbb{Z}$-graded Lie algebra.
If $\frak g=\oplus_{-r\leq i\leq s}\frak g_{_{[i]}}$ is a $\mathbb
    Z$-graded Lie  algebra, then $\oplus_{-r\leq i\leq 0}\frak
    g_{_{[i]}}$ is called the \textit{top of} $\frak g$ (with respect to the
    gradation). Let $V=\oplus_{i\in \mathbb{Z}}V_{[i]}$ be a $\mathbb{Z}$-graded vector space and $x\in V$ is
    a $\mathbb{Z}$-homogeneous element. Then we let $\mathrm{zd}(x)$ denote the $\mathbb{Z}$-degree of $x.$

In the following we recall the notions  of the generalized Witt
modular superalgebra and the Hamiltonian modular superalgebra and
their natural  gradation structures. We  also  introduce the notation,
terminology and convention which will be used throughout.

In the sequel  $\mathbb {F}$ denotes a field of characteristic
$p>3.$ In addition to the standard notation $\mathbb {Z},$ we use
 $\mathbb {N}$ for the set of positive integers and
${\mathbb {N}}_0$ for the set of nonnegative integers. Henceforth,
we will let  $m$ and $n$  denote fixed positive integers  without notice. Given $\alpha = (\alpha _1,\ldots,\alpha
_m ) \in \mathbb {N}_0^m,$ we put $|\alpha| :=\sum_{i=1}^m\alpha
_i.$ Following \cite{st}, denote by $\mathcal{O}(m)$ the \textit{divided power algebra
over}
 $\mathbb{F}$ with an ${\mathbb
F}$-basis $\{ x^{( \alpha ) }\mid \alpha \in \mathbb{N}_0^m \}.$
For $\varepsilon _i=( \delta _{i1},\ldots,\delta _{im}),$ we
abbreviate $x^{( \varepsilon _i)}$ to $x_i,$ $i=1,\ldots,m.$ Let
$\Lambda (n)$ be the \textit{exterior superalgebra
over} $\mathbb{F}$ in $n$ variables
$x_{m+1},\ldots, x_{m+n}.$   Denote the tensor product by
$\mathcal{O}(m,n) =\mathcal{O}(m)\otimes_{\mathbb{F}} \Lambda(n).$
Obviously, $\mathcal{O}(m,n)$ is an associative superalgebra with a
$\mathbb{Z}_2$-gradation induced by the trivial
$\mathbb{Z}_2$-gradation of $\mathcal{O}(m)$ and the natural
$\mathbb{Z}_2$-gradation of $\Lambda (n).$ Evidently,
$\mathcal{O}(m,n)$ is super-commutative.

For $g\in \mathcal{O}(m),f\in \Lambda(n),$ we abbreviate $ g\otimes
f$ to $gf.$ The following formulas hold in $\mathcal{O}(m,n):$
$$
x^{(\alpha) }x^{(\beta) }=\binom{\alpha +\beta }{ \alpha}
 x^{( \alpha +\beta)}\quad\mbox{for}\  \alpha,\beta \in
{\mathbb N}_0^m;
$$

$$
x_ix_j=-x_jx_i\quad\mbox{for}\ i,j=m+1,\ldots, m+n;
$$
$$
x^{( \alpha ) }x_j=x_jx^{( \alpha ) }\quad\mbox{for}\ \alpha \in
\mathbb{N}_0^m,\, j=m+1,\ldots, m+n,
$$
where $\binom{ \alpha +\beta} {\alpha}:=\prod_{i=1}^m\binom{
\alpha _i+\beta _i}{ \alpha _i}.$\vspace{2mm}

 For convenience,
put $Y_0:=\{ 1,2,\ldots,m \},$ $Y_1:=\left\{m+1,\ldots,
m+n\right\} $ and $Y:=Y_0\cup Y_1.$ Let
$$\mathbb{B}_k:=\left\{ \langle
i_1,i_2,\ldots,i_k\rangle |m+1\leq i_1<i_2<\cdots <i_k\leq
m+n\right\}
$$
be the set of $k$-tuples of strictly increasing integers between $m+1$ and $m+n,$ and $\mathbb{B}
:=\bigcup\limits_{k=0}^n\mathbb{B}_k,$ where
$\mathbb{B}_0:=\emptyset.$ For $u:=\langle i_1,i_2,\ldots
,i_k\rangle \in \mathbb{B}_k,$ set $ |u| :=k,| \emptyset | :=0$,
$x^\emptyset :=1,$ and $x^u:=x_{i_1}x_{i_2}\cdots x_{i_k};$ we
also let $u$ stand for the set $\{i_1,i_2,\ldots,i_k\} $ if no
confusion occurs.
Put $\mathbb{B}^{0}:=\{u\in \mathbb{B}\mid |u|\, \mbox{even}\}.$
 Clearly, $\left\{ x^{\left( \alpha \right)
}x^u\mid\alpha \in \mathbb{N}_0^m,u\in \mathbb{B}
\right\} $ is an $\mathbb{F}$-basis of $\mathcal{O} \left(
m,n\right).$

Let $\partial_1,\partial_2,\ldots,\partial_{m+n}$ be the linear transformations of
$\mathcal{O} \left( m,n\right) $ such that
\[
\partial_i ( x^{( \alpha ) }x^u) =\left\{
\begin{array}{l}
x^{\left( \alpha -\varepsilon _i\right) }x^u,\quad \quad \quad
\quad i\in Y_0
\\
x^{\left( \alpha \right) }\, \partial x^u/\partial x_i,\quad
\quad i\in Y_{1.}
\end{array}
\right.
\]
Then $\partial_1,\partial_2,\ldots,\partial_{m+n}$ are superderivations of the
superalgebra $\mathcal{O} \left( m,n\right).$ Let
\[
W\left(m,n\right) :=\Big\{ \sum\limits_{i\in Y}a_i \partial_i \mid a_i\in
\mathcal{O} \left( m,n\right),i\in Y\Big\}.
\]
Then $W\left( m,n\right) $ is a Lie superalgebra contained in
${\rm Der}\,  \mathcal{O} \left(m,n\right). $

One may verify that
\begin{equation}[aD,bE]=aD(b)E-(-1)^{{\rm
p}(aD){\rm p}(bE)}bE(a)D+(-1)^{{\rm p}(D){\rm p}(b)}ab[D,E]\label{le0}
\end{equation}
for $ a,b\in \mathcal{O} (m,n),$ $ D,E\in {\rm Der} \mathcal{O} (m,n
).$ Consequently, the following formula holds in $W\left(
m,n\right):$
\begin{equation*}
\left[ a\partial_i,b\partial_j\right] =a\partial_i\left( b\right) \partial_j-\left( -1\right)
^{{\rm p}(a\partial_i) {\rm{p}}(b\partial_j) }b\partial_j\left( a\right) \partial_i
\end{equation*}
for $a,b\in \mathcal{O} \left( m,n\right),i,j\in Y.$ We note that
${\rm p}( \partial_i) =\mu (i),$ where
\[
\mu \left( i\right) :=\left\{
\begin{array}{l}
\overline{0},\quad \quad i\in Y_0 \\
\overline{1},\quad \quad i\in Y_{1.}
\end{array}
\right.
\]

In the sequel suppose    $m=2r$ is even. Define the linear mapping $
\mathrm{D_{H}}:\mathcal{O} ( m,n)\rightarrow W(m,n)$ by means of
                $$ \mathrm{D_{H}}(a):=\sum _{i\in Y} \tau(i)(-1)
               ^{\mu(i) {\mathrm{p}}(a)} \partial_{i}(a) \partial_{i'}
                 \quad\mbox{for all}\; a \in \mathcal{O} (m,n), $$
 where $$
i':=\left\{
\begin{array}{l}
      i+r,\quad 1\leq i \leq r\\
      i-r,\quad r< i \leq 2r\\
      i, \quad \quad \quad i \in Y_{1};
      \end{array} \right.
      \quad\quad \quad \tau(i):=\left\{ \begin{array}{l}
      \ \ 1,\quad 1\leq i \leq r\\
      -1,\quad r< i \leq 2r\\
      \ \ 1,\quad i \in Y_{1}.
      \end{array}\right.
     $$Then the following identity holds
  \begin{equation} \label{le1} \big [\mathrm{D_{H}}(a),\mathrm{D_{H}}(b)\big ]=\mathrm{D_{H}}
  (\mathrm{D_{H}}(a)(b))\quad
  \mbox{for all}\     a, b \in \mathcal{O} (2r,n).
  \end{equation}

Let $$ \underline{t}:=( t_{1},t_{2},\ldots,t_m) \in \mathbb{N}^m,\
\  \pi :=( \pi _1,\pi _2,\ldots,\pi _m)$$ where
 $\pi _i:=p^{t_i}-1,i\in Y_0.$ Put $\mathbb{A}:=\mathbb{A}(m;\underline{t})
:=\left\{ \alpha \in \mathbb{N}_0^m|\alpha _i\leq \pi
_i,i=1,2,\ldots,m\right\}.$ Then
$$ \mathcal{O} \left( m,n;\underline{t}\right) :={\mathrm
{span}}_{\mathbb{F}}\big\{ x^{\left( \alpha \right) }x^u \mid
\alpha \in \mathbb{A},u\in \mathbb{B} \big\}
$$ is a finite-dimensional subalgebra of $\mathcal{O}
\left(m,n\right),$ with a natural ${\mathbb Z}$-gradation
$\mathcal{O} \left(m,n;\underline{t}\right)=\oplus _{i=0}^{\xi}
\mathcal{O}(m,n;\underline{t})_{[i]},$ where
        $\mathcal{O}(m,n;\underline{t})_{[i]}=
        {\mathrm {span}} _{\mathbb F}\{ x^{(\alpha)} x^{u}\mid |\alpha|+|u|=i \},
         \xi:=|\pi|+n.$
Set
$$ W\left( m,n;\underline{t}\right) :=\Big\{
\sum_{i\in Y}a_i\partial_i \mid a_i\in \mathcal{O} \left(
m,n;\underline{t}\right),i\in Y\Big\}.
$$
Then $W\left( m,n;\underline{t}\right) $ is a  subalgebra of
$W\left( m,n\right).$ In particular, it is a finite-dimensional
simple Lie superalgebra (see \cite{Zh1}). As in the case of Lie
algebras, $W\left( m,n;\underline{t}\right) $ is called
\textit{the generalized Witt superalgebras}. Obviously, $W\left(
m,n;\underline{t}\right) $ is a free $ \mathcal{O} \left(
m,n;\underline{t}\right) $-module with $\mathcal{O} \left(
m,n;\underline{t}\right) $-basis $\left\{ \partial_1,\partial_2,\ldots
,\partial_{m+n}\right\}.$

Set
$$
H(m,n;\underline{t}):=\big\{ \mathrm{D_H}\left( a\right) \mid a\in
\oplus_{i=0}^{\xi-1}\mathcal{O} \left( m,n;\underline{t}\right)_{[i]}
\big\},$$ where $\xi=|\pi|+n.$ Then $H( m,n;\underline{t})$ is a
finite-dimensional simple Lie superalgebra,
 which is called  \textit{the Hamiltonian
 superalgebra} (see \cite{ZF,Zh1}).

The $\mathbb{Z}$-gradation of    $\mathcal{O}(m,n;\underline{t})$
induces naturally a  $\mathbb{Z}$-gradation structure   of the
generalized Witt  superalgebra $W(m,n;\underline{t})
=\oplus_{i=-1}^{\xi-1} W(m,n;\underline{t})_{[i]},$  where
$$W(m,n;\underline{t})_{[i]} :={\rm span}_{\mathbb{F}}\{f\partial_{s}\mid
s\in Y,\ f\in \mathcal{O}(m,n;\underline{t})_{[i+1]}\}.$$
Note that   $H(m,n;\underline{t})$ is a  $\mathbb{Z}$-graded subalgebras of
$W(m,n;\underline{t}).$

 In the following sections, since the positive integer $m$ is even, we write $2m$ instead of $m$,
 and we usually write $\mathcal{O},$ $W$,
and $H$ for $\mathcal{O}(2m,n;\underline{t}),$
 $W(2m,n;\underline{t}),$
 and $H(2m,n;\underline{t})$, respectively. For convenience, the even parts of
 $W$
and $H$ will be denoted by $ \mathcal{W}$ and $\mathcal{H}$, respectively.

In this paper we suppose   $m>1$, $n>3$ and $p>3$ for the sake of
simplicity although sometimes a weak hypothesis is sufficient.

\section{The ideal $\frak{N}$ and its generating set }

As mentioned in the introduction, our main object is to discuss the
derivations for the even part $\mathcal{H}$ of the finite-dimensional
Hamiltonian  superalgebra $H.$ Precisely speaking, we want to
formulate certain homogeneous derivations from $\mathcal{H}$ into
$\mathcal{W}.$ However, in contrast to the setting of
$\mathcal{W}$ and $\mathcal{S} $ (see \cite{lg,lz1}), we shall first deal with the
derivations from an ideal $\frak{N}$  of $\mathcal{H}$ into $\mathcal{W}$ rather
than the derivations  from $\mathcal{H}$ into $\mathcal{W}.$ Of course, we can
guarantee that this transformation  behaves without any  influence
on our aim (see Remark \ref{r1.5}).

Recall
$$\mathcal{H}:=H_{\overline{0}}(2m,n;\underline{t})={\rm
span}_{\mathbb{F}}\big\{  \mathrm{D_{H}}(x^{(\alpha)}x^u)\mid
\alpha\in \mathbb{A},u\in
\mathbb{B}^0,(\alpha,u)\not=(\pi,\omega)\big\}
$$
and put
$$\frak{N}:={\rm span}_{\mathbb{F}}\big\{
\mathrm{D_{H}}(x^{(\alpha)}x^u)\ |\ \alpha\in \mathbb{A},u\in
\mathbb{B}^0
,(\alpha,u)\not=(\pi,\omega),(\alpha,u)\not=(\pi,\emptyset)\big\},
$$
where $\mathbb{B}^0=\{u\in \mathbb{B}\mid |u|\  \mbox{even}\} $ and
$\omega:=\langle m+1,\ldots,m+n\rangle\in \mathbb{B}_{n}.$
Evidently,  $\frak{N}$ is a subspace of
 $\mathcal{H}$ of codimension 1:
\begin{equation}
\mathcal{H}=\frak{N}\oplus\mathbb{F}\,
\mathrm{D_{H}}(x^{(\pi)}).\label{e6.6.1}
\end{equation}
In the following we shall demonstrate that in order to determine the
derivations from $\mathcal{H}$ into $\mathcal{W} $ it suffices to
determine the derivations from $\frak{N}$ into $\mathcal{W}.$ It is
clear that this will  simplify our consideration. To that end, we
need the following two propositions.

\begin{proposition}\label{p1.1}
$\frak{N}$ is an ideal of $\mathcal{H}.$
\end{proposition}

 \begin{proof}  We first show that $\frak{N}$ is a
subalgebra of $\mathcal{H}.$ Given any two linear generators
$\mathrm{D_{H}}(x^{(\alpha)}x^u) $,
$\mathrm{D_{H}}(x^{(\beta)}x^v)\in \frak{N},$ we assert that
\begin{equation}
[\mathrm{D_{H}}(x^{(\alpha)}x^u), \mathrm{D_{H}}(x^{(\beta)}x^v)
]\in \frak{N}.\label{e6.6.2}
\end{equation}
In fact, if $u=v=\emptyset,$ in view of the theory of Hamiltonian
algebras \cite{SF}, we get (\ref{e6.6.2}) immediately. If
  $u\not=\emptyset$ or
  $v\not=\emptyset,$ noticing that
   $|u|,$ $|v|$ are all even, we obtain (\ref{e6.6.2}) by using the formula
   (\ref{le1}).
%$$[\mathrm{D_{H}}(f),\mathrm{D_{H}}(g)]=\mathrm{D_{H}}(\mathrm{D_{H}}(f)(g)),$$
 This proves that
  $\frak{N}$ is a subalgebra of
   $\mathcal{H}. $ We next show that
    $\frak{N}$ is an ideal of
    $\mathcal{H}.$ By (\ref{e6.6.1}), it suffices to show that
\begin{equation}
[\mathrm{D_{H}}(x^{(\pi) }),\frak{N} ]\subset
\frak{N}.\label{e6.6.3}
\end{equation}
Clearly,
\begin{equation} [\mathrm{D_{H}}(x^{(\pi) }),\frak{N}
]\subset\mathcal{H}.\label{e6.6.4}
\end{equation}
 For every linear generator
$\mathrm{D_{H}}(x^{(\alpha)}x^u)$ of $\frak{N},$ if
 ${\rm zd} ([\mathrm{D_{H}}(x^{(\pi) }), \mathrm{D_{H}}(x^{(\alpha)}x^u)])\not= |\pi|-2,$
then we see from (\ref{e6.6.4}) th at $[\mathrm{D_{H}}(x^{(\pi) }),
\mathrm{D_{H}}(x^{(\alpha)}x^u)]\in \frak{N} $ and thereby
(\ref{e6.6.3}) holds. It remains only to consider the case
$|\alpha|+|u|=2.$ If
 $|u|=2,$ then
   $[\mathrm{D_{H}}(x^{(\pi) }), \mathrm{D_{H}}(x^u)]=0\in
\frak{N};$ if
  $\alpha=2,$   by the theory of Hamiltonian
  Lie algebras, we also have
$[\mathrm{D_{H}}(x^{(\pi) }), \mathrm{D_{H}}(x^{(\alpha)})]\in
\frak{N} $ and therefore, (\ref{e6.6.3}) holds.
\end{proof}

Following \cite{lz1}, put $\mathcal{G}:=\mathrm{span}_\mathbb{F}\{x^u\partial_r\,|\,r\in Y,u\in
\mathbb{B},\mathrm{p}(x^u\partial_r)=\overline{0}\}.$
Then $C_{\mathcal{W} } (\mathcal{W} _{-1})=\mathcal{G};$ in particular,
$C_{\mathcal{H} } (\mathcal{H} _{-1})\subset\mathcal{G}.$ Note that $\mathcal{G}$
is a $\mathbb{Z}$-graded subalgebra of $\mathcal{W}. $

\begin{proposition}\label{p1.2} Let $|\omega|=n$ be even. Then
$C_{\mathcal{W}}(\frak{N})=\mathbb{F}\,
\mathrm{D_{H}}(x^{\omega}).$
\end{proposition}
 \begin{proof}
 For any arbitrary basis element
 $\mathrm{D_{H}}(x^{(\alpha)}x^u)$ of $\frak{N},$
 $u\in \mathbb{B}^0,$ it is clear that
 $$[\mathrm{D_{H}}(x^{\omega}),\mathrm{D_{H}}(x^{(\alpha)}x^u)]=0.$$
 Therefore, $\mathbb{F}\,
\mathrm{D_{H}}(x^{\omega})\subset
C_{{\mathcal{W}}}(\frak{N}). $
 We propose to prove the converse inclusion.
 Since $C_{ \mathcal{W} }(\frak{N})$ is a
 $\mathbb{Z}$-subalgebra of $ \mathcal{H},$
 it suffices to show that if
 $D\in C_{\mathcal{W}}(\frak{N})$ is homogeneous  then
 $D\in \mathbb{F}\,
\mathrm{D_{H}}(x^{\omega}).$ Noting that
 $\mathcal{W}_{[-1]}=\mathcal{H}_{[-1]}=\frak{N}_{[-1]},$ we see
 that
 $D\in \mathcal{G}.$ Thus one may assume that
\begin{equation}
D=\sum_{r\in Y} f_r\partial_r\quad {\rm where}\ f_r \in
\Lambda(n).\label{e6.6.17}
 \end{equation}
 For any $i\in Y_0,$ since
 $\mathrm{D_{H}}(x^{(2\varepsilon_i)})\in \frak{N},$
 we have
 $[D,\mathrm{D_{H}}(x^{(2\varepsilon_i)})]=0 $
 and therefore,
 $f_i\partial_{i'}=0.$
 This proves that $f_i=0$ for all
 $i\in Y_0.$
 Thus, by (\ref{e6.6.17}), we have
\begin{equation}
D=\sum_{r\in Y_1} f_r\partial_r\quad\mbox{where}\ f_{r}\in
\Lambda(n).\label{e6.6.18}
 \end{equation}

   \noindent\textit{Case} (i): $1<\mathrm{zd}(f_{r})\leq n-3 $ for
 $r\in Y_{1}.$
Assume that $f_{r_{_0}}\not= 0$ for some $r_{_0}\in Y_1.$
 By our assumption   one may choose
    $k,l\in Y_1$ with
  $k\not=l$ such that
  $x_l\partial_k(f_{r_{_0}})\not=0 $ and $k\neq r_{_0},$  $l\neq r_{_0}.$
  Note that $\mathrm{D_{H}}(x_kx_l)=x_l\partial_k-x_k \partial_l\in\frak{N}.$
  We have
 \begin{equation}
0=\Big[\mathrm{D_{H}}(x_kx_l),\sum_{r\in Y_1} f_r \partial_r\Big]
=\sum_{r\in Y_1} (x_l \partial_k(f_r)-x_k \partial_l(f_r)
)\partial_r-f_l \partial_k+f_k \partial_l.\label{e6.6.19}
 \end{equation}
    Then a comparison of the coefficients of $\partial_{r_{_0}}$ in (\ref{e6.6.19}) yields that
   $x_l\partial_k(f_{r_{_0}})-x_k\partial_l(f_{r_{_0}})=0.$ However, this equation
   implies
   $x_l\partial_k(f_{r_{_0}})=0, $ a contradiction. This
   proves that
   $D=0.$\newline

    \noindent\textit{Case} (ii): $\mathrm{zd}(f_{r})=n-1$ for
 $r\in Y_{1}.$ Then one may assume that
  \begin{equation}\label{ee6.6.19}
   f_{r}=\sum_{s\in Y_1}c_{rs}x^{\langle
   \widehat{s}\rangle}\quad\mbox{where} \ c_{rs}\in \mathbb{F}.
  \end{equation}
   Here we put $\langle
   \widehat{s}\rangle:=\omega-\langle
    s \rangle.$ Fixing $r\in Y_1$ and replacing $f_{r}$ in (\ref{e6.6.19})  by
    (\ref{ee6.6.19}),  for $k, l\in Y_{1}\setminus r$
    with $k\neq l,$ we obtain that $x_{l}\partial_{k}(\sum_{s\in Y_1}c_{rs}x^{\langle
   \widehat{s}\rangle})=0$ and therefore, $c_{rl}x^{\langle
   \widehat{k}\rangle}=0.$ This implies that  $c_{rl}=0$  whenever $r\neq
   l.$ Consequently, $f_{r}=c_{rr}x^{\langle
   \widehat{r}\rangle} $ for $r\in Y_1.$ Observing the
   coefficient of $\partial_{k}$ in (\ref{e6.6.19}), we obtain that
   $-x_k\partial_l(c_{kk}x^{\langle
   \widehat{k}\rangle})-c_{ll}x^{\langle
   \widehat{l}\rangle}=0. $ It follows that $c_{kk}=(-1)^{k+l}c_{ll}
   $ for $k,l\in Y_{1}.$ Let $\lambda:=c_{m+1,m+1}.$ So far, we
   have proved that
   $$D=\lambda \sum_{r\in Y_1}(-1)^{r+1}x^{\langle \widehat{r}\rangle}
   =\lambda \mathrm{D_H}(x^{\omega}).$$

\noindent\textit{Case} (iii): $\mathrm{zd}(f_{r})=1.$ As in Case (i) one may easily get the desired result.

   \end{proof}

   \begin{proposition}\label{p1.3} Let $|\omega|=n$ be odd. Then
   $C_{\mathcal{W}}(\frak{N})=0.$
   \end{proposition}
    \begin{proof} Arguing just as in the proof of
   Proposition \ref{p1.2}
    we see that every element $D\in
   C_{\mathcal{W}}(\frak{N})$ may be written  as (\ref{e6.6.18}).
   Given $r_{_{0}}\in Y_{1},$ we obtain from (\ref{e6.6.19})
   that
   \begin{equation}x_{l}\partial_{k}(f_{r_{_{0}}})=0 \quad
   \mbox{whenever} \,\, k, l\in
   Y_{1}\setminus r_0\,\, \mbox{with}\, \, k\neq l.\label{0510}
   \end{equation}
   Without loss of generality, one may
   assume that $D $ is homogeneous and so is $f_{r_{_{0}}}=\sum_{u}
   c_{r_{_{0}}u}x^{u}.$ If $1<\mathrm{zd}(f_{r_{_{0}}})<n-1, $  then by (\ref{0510}),
    $f_{r_{0}}=0.$ In the case $\mathrm{zd}(f_{r_{_{0}}})=1,$ it is easily shown that $D=0.$
    Assume that $f_{r_{_{0}}}\neq 0.$
   Then $\mathrm{zd}(f_{r_{_{0}}})\geq n-1 $ as we have shown that
   $x_{l}\partial_{k}(f_{r_{_{0}}})=0 $
   whenever $k, l\in
   Y_{1}$ with $k\neq l.$ Because $n$ is odd, this forces
   $\mathrm{zd}(f_{r_{_{0}} })=n $  and hence, $f_{r_{_{0}}}=
   c_{r_{_{0}}\omega}x^{\omega}.$ Now, using (\ref{e6.6.19}) one may
   easily deduce that $c_{r_{_{0}}\omega}=0,$ contradicting the
   assumption that $f_{r_{_{0}}}\neq 0.$ Summarizing, we have
   shown that $D=0.$ The proof is complete.\end{proof}

   \begin{theorem}\label{t1.4} Suppose
   $\varphi\in \mathrm{Der}(\mathcal{H}, \mathcal{W}) $ and
   $\varphi(\frak{N})=0.$ Then the following statements hold.

    $\mathrm{(i)}$\; If $|\omega|=n$ is odd then $\varphi=0.$

   $\mathrm{(ii)}$\ If $|\omega|=n$ is even then
    $\varphi
   (\mathrm{D_H}(x^{(\pi)}))=\lambda \mathrm{D_H}(x^{\omega}) $
   for some $\lambda\in \mathbb{F}.$
   Conversely, any mapping
   $\varphi:\mathcal{H}\rightarrow\mathcal{W}$ vanishing on
   $\frak{N}$ and satisfying  $\varphi
   (\mathrm{D_H}(x^{(\pi)}))=\lambda \mathrm{D_H}(x^{\omega})$ for
   any fixed $\lambda\in \mathbb{F}$ is necessarily a derivation
   from $\mathcal{H}$ into $\mathcal{W}.$
   \end{theorem}
 \begin{proof} By Proposition
     \ref{p1.1},
    $\big[\varphi(\mathrm{D_H}(x^{(\pi)})),\frak{N}\big]=0;$
    that is, $\varphi(\mathrm{D_H}(x^{(\pi)}))\in
    C_{\mathcal{W}}(\frak{N}).$ Now (i) follows from Proposition
    \ref{p1.3}
    and the decomposition (\ref{e6.6.1}). The first part of
    (ii) is an immediate consequence of Propositions \ref{p1.1}, \ref{p1.2}  and
    the decomposition (\ref{e6.6.1}). The second part follows
    immediately from Propositions \ref{p1.1} and \ref{p1.2}.
    \end{proof}

    \begin{remark}\label{r1.5} Suppose that the structure of $\mathrm{Der}(\frak{N},\mathcal{W})$ has been determined and
    all the derivations from $\frak{N}$ into $\mathcal{W}$ may extend to $\mathcal{H}.$
    Then, in the light of Theorem \ref{t1.4},  one may easily
    determine the derivation space $\mathrm{Der}(\mathcal{H},\mathcal{W}).$
    We thereby pay our attention to the ideal $\frak{N}$ in
    future.
    \end{remark}

   \begin{remark}\label{r1.6}  Define the linear mapping
   $\Gamma_{\lambda}:\mathcal{H}\rightarrow\mathcal{W}$  by means of
   $\Gamma_{\lambda}(\frak{N}):=0$  and $\Gamma_{\lambda}
   (\mathrm{D_H}(x^{(\pi)})):=\lambda \mathrm{D_H}(x^{\omega})$ for
   any fixed $\lambda\in \mathbb{F}.$ By Theorem \ref{t1.4},
   $\Gamma_{\lambda}$ is a derivation from $\mathcal{H}$ into $\mathcal{W}.$
   Note that $\mathrm{zd}(\Gamma_{\lambda})=|\omega|-|\pi|.$ In general,
   $\Gamma_{\lambda}$ is outer. In the next section, we shall give in addition three series
   of the so-called exceptional derivations from $\mathcal{H}$ into $\mathcal{W} $
   in the setting that $|\omega|=n$
   is even. In general, they are outer derivations.
  \end{remark}

In the following we study the  generating  set of $\frak{N}.$ Put
$$\mathcal{M}:=\{ \mathrm{D_{H}}(x^{(q_i\varepsilon_i)})\mid
1\leq q_i \leq \pi_i,\ i\in Y_0\},$$ and
$$\mathcal{N}:=\{ \mathrm{D_{H}}(x_ix^u) \mid i\in Y_0,u\in \mathbb{B}_2 \}.$$
We conclude this section with the following generating theorem.

\begin{theorem}\label{t1.7}
$\frak{N}$ is generated by $\mathcal{M}\cup \mathcal{N} \cup
\frak{N}_{[0]}.$
\end{theorem}
  \begin{proof} Let $\mathcal{L}$ be the subalgebra of
$\frak{N} $ generated by $\mathcal{M}\cup \mathcal{N} \cup
\frak{N}_{[0]}.$ For $i,j\in Y_0,$ direct computation shows that
\begin{equation}
\mathrm{D_{H}}(x^{(2\varepsilon_i)
}x_{j'})=\tau(i)(1+\delta_{ij})^{-1}[\mathrm{D_{H}}(x^{(3\varepsilon_i)
}),\mathrm{D_{H}}(x_{i'} x_{j'}) ]\in \mathcal{ L}.\label{e6.6.5}
\end{equation}
Furthermore,
\begin{equation}
\mathrm{D_{H}}(x_ix_{i'}x_{j'})=\tau(i')(1-\frac{1}{3}\delta_{ij'})
[\mathrm{D_{H}}(x^{(2\varepsilon_{i'})}),
\mathrm{D_{H}}(x^{(2\varepsilon_{i})} x_{j'}) ]\in \mathcal{
L}.\label{e6.6.6}
\end{equation}
Using (\ref{e6.6.6}) we get for $i\not=j',$
%$i\not=j'$
\begin{equation}
\mathrm{D_{H}}(x^{(\pi_i\varepsilon_i)}
x_{j'})=-\tau(i)(1+\delta_{ij})^{-1}[\mathrm{D_{H}}(x^{(\pi_i\varepsilon_{i})}),
\mathrm{D_{H}}(x_{i}x_{i'} x_{j'}) ]\in \mathcal{
L}.\label{e6.6.7}
\end{equation}
It follows from (\ref{e6.6.5}) and (\ref{e6.6.7}) that for
$i\not=j,j',$
\begin{equation}
\mathrm{D_{H}}(x^{(\pi_i\varepsilon_i)}x_{j}
x_{j'})=\tau(j)[\mathrm{D_{H}}(x^{(2\varepsilon_{j})}x_{j'}),
\mathrm{D_{H}}(x^{(\pi_i \varepsilon_i)} x_{j'}) ]\in \mathcal{
L}.\label{e6.6.8}
\end{equation}
An application of (\ref{e6.6.8}) yields
$$\mathrm{D_{H}}(x^{(\pi_1\varepsilon_1+\pi_2\varepsilon_2)}
)=-[\mathrm{D_{H}}(x^{(\pi_1\varepsilon_1)}),
\mathrm{D_{H}}(x_1x_{1'}x^{(\pi_2\varepsilon_2)})]\in \mathcal{
L};$$ and therefore,
$$\mathrm{D_{H}}(x^{(\pi_1\varepsilon_1+\pi_2\varepsilon_2+\pi_3\varepsilon_3)}
)=-[\mathrm{D_{H}}(x^{(\pi_1\varepsilon_1+\pi_2\varepsilon_2)}),
\mathrm{D_{H}}(x_1x_{1'}x^{(\pi_3\varepsilon_3)}) ]\in \mathcal{
L}.$$ By induction we may easily obtain that
\begin{equation}
\mathrm{D_{H}}(x^{(\pi_1\varepsilon_1+\cdots+\pi_m\varepsilon_m)}
)\in \mathcal{ L}.\label{e6.6.9}
\end{equation}
Using (\ref{e6.6.9}) and (\ref{e6.6.8}) we obtain that
$$\mathrm{D_{H}}(x^{(\pi_1\varepsilon_1+\cdots+\pi_m\varepsilon_m+\pi_{1'}
\varepsilon_{1'})}
)=-[\mathrm{D_{H}}(x^{(\pi_1\varepsilon_1+\cdots+\pi_m\varepsilon_m)}
),\mathrm{D_{H}}(x_mx_{m'}x^{(\pi_1'\varepsilon_{1'})}) ]\in
\mathcal{ L}.$$ By induction one may easily show that
\begin{equation}
\mathrm{D_{H}}(x^{(\pi-\pi_{2m} \varepsilon_{2m})})\in\mathcal{
L}.\label{e6.6.10}
\end{equation}
By (\ref{e6.6.7}),  $\mathrm{D_{H}}(x_mx^{(\pi_{m'}
\varepsilon_{m'})})\in\mathcal{ L}.$ It follows from
(\ref{e6.6.10}) that
\begin{equation}
\mathrm{D_{H}}(x^{(\pi-\varepsilon_{m'}
)})=-[\mathrm{D_{H}}(x^{(\pi-\pi_{m'}\varepsilon_{m'} )}),
\mathrm{D_{H}}(x_m x^{(\pi_{m'}\varepsilon_{m'} )})] \in\mathcal{
L}.\label{e6.6.11}
\end{equation}
In general, we have
$$\mathrm{D_{H}}(x^{(\pi-\varepsilon_{r})})\in\mathcal{ L}
\quad  \mbox{for all} \  r\in Y_0. $$ Therefore,
\begin{equation}
\mathrm{D_{H}}(x^{(\alpha)})\in\mathcal{ L} \quad \mbox{for all}\
\alpha\not=\pi.\label{e6.6.12}
\end{equation}
We want to prove that
\begin{equation}
\mathrm{D_{H}}(x^{(\alpha)}x^u)\in\mathcal{ L} \quad \mbox{for}\
\alpha\not=\pi,\ {\emptyset}\not=u\in
\mathbb{B}^{0}.\label{e6.6.13}
\end{equation}
We proceed by induction on $|\alpha|+|u|.$ When $|\alpha|+|u|=2,$
we have
 $u\in \mathbb{B}_2$ and $\alpha=0 $ and therefore,
$\mathrm{D_{H}}(x^u)\in\mathcal{ L},$ since
$\mathcal{H}_{[-1]}\cup \mathcal{N}\subset \mathcal{ L}; $ that
is, (\ref{e6.6.13}) holds. When $|\alpha|+|u|=3,$ necessarily
$|u|=2$ and $|\alpha|=1$ and therefore, (\ref{e6.6.13}) holds.
Assume that $|\alpha|+|u|>3. $
  Find $w \in\mathbb{B}_2$ and $v\in \mathbb{B}$ such that $x^{v}x^{w}=x^{u}.$
  If  $\alpha_i<\pi_i $ for some $i\in Y_{0},$ then
 $|\alpha+\varepsilon_i|+|v|<|\alpha|+|u|,$
 $1+|w|<|\alpha|+|u|.$
 Thus,
 $$\mathrm{D_{H}}(x^{(\alpha)}x^u)
 =\tau(i)[\mathrm{D_{H}}(x^{(\alpha+\varepsilon_i)}x^v),\mathrm{D_{H}}(x_{i'}x^{w})]
 \in\mathcal{ L};$$
 that is, (\ref{e6.6.13}) holds.

 It remains to show that
\begin{equation}
\mathrm{D_{H}}(x^{(\pi)}x^u)\in\mathcal{ L} \quad {\rm for} \ u\in
\mathbb{B}\  \mbox{ with} \ u\not= \omega, u
\not=\emptyset.\label{e6.6.14}
\end{equation}
To do that, we first  show that
\begin{equation}\mathrm{D_{H}}(x^u)\in \mathcal{L}\quad \mbox{for all}\ u\in
\mathbb{B}^0.\label{e6.6.15}
\end{equation}
Clearly, $\mathrm{D_{H}}(x^u)\in\mathcal{ L} $ for all $u\in
\mathbb{B}_2$. Assume that
  $u\in \mathbb{B}^0$ and $|u|>2.$
 Find  $v, w \in\mathbb{B}^0 $ satisfying $|v|<|u| $ and $|w|<|u|$, such that $x^{v}x^{w}=x^{u}.$
  Then
  \begin{equation*}\mathrm{D_{H}}(x^u ) =
  [\mathrm{D_{H}}(x_1x^v),\ \mathrm{D_{H}}(x_{1'}x^w)].
  \end{equation*}
  Since $\mathcal{N}\subset \mathcal{L},$ using
  induction one may easily prove (\ref{e6.6.15}).
 When $|u|=2,$ let $u=\langle k,l\rangle.$ Take $r\in Y_1\setminus\{k,l\}.$
 By (\ref{e6.6.13}),
\begin{equation}
\mathrm{D_{H}}(x^{(\pi)}x_kx_l)
=[\mathrm{D_{H}}(x^{(\pi-\pi_1\varepsilon_1)}x_rx_k),\
\mathrm{D_{H}}(x^{(\pi_1\varepsilon_1)}x_rx_l) ]
 \in\mathcal{ L};\label{e6.6.16}
\end{equation}
that is,  (\ref{e6.6.14}) holds in this case. Now suppose
 $|u|>2$ and $u\not=\omega.$ Find $v,w\in \mathbb{B}$ with $ |v|=1  $ such that $x^{v}x^{w}=x^{u}.$
 Take
 $r\in Y_1\setminus u.$ Then (\ref{e6.6.15}) and (\ref{e6.6.16}) ensure that
 $$\mathrm{D_{H}}(x^{(\pi)}x^u)
 =[\mathrm{D_{H}}(x^{(\pi)}x_rx^v), \mathrm{D_{H}}( x_rx^w) ]\in \mathcal{ L}, $$
 proving (\ref{e6.6.14}). The proof is complete.
 \end{proof}

\section{Exceptional derivations}

  In the this section we shall give three series of the so-called \textit{exceptional derivations
from}  $\mathcal{H}$ into $\mathcal{W}.$ As
we shall see, in general these exceptional derivations are all
outer. We note that these exceptional derivations have no analogs
in the setting   for the even part  of the
odd Hamiltonian superalgebra or the special  superalgebra (see \cite{lhs,lz1}). Thus,
roughly speaking,   the even part of Hamiltonian  modular Lie
superalgebra possesses ``more" outer derivations than the even part of
 the special Lie superalgebra. The other reason we are
interested in this phenomenon is that it does not occur in the
``super" setting (see \cite[Theorem 2.13]{WZ3}).   Throughout this section assume
that $|\omega|=n$ is even.

Let us define the first series of exceptional derivations from $\mathcal{H}$
into $\mathcal{W}$.  Given $q\in \mathbb{N}$ and $i\in Y_0,
$ define
\begin{eqnarray*}
\Phi^{(q)}_{i}: \mathcal{H} \rightarrow \mathcal{W},\quad
 \mathrm{D_{H}}(f)\mapsto \partial_i^{p^q}(f)\mathrm{D_{H}}(x^{\omega}).
 \end{eqnarray*}
 As $\mathrm{Ker}(\mathrm{D_{H}})=\mathbb{F}\, 1,$ $\Phi^{(q)}_{i}$
 is well defined. Clearly,
${\rm zd}(\Phi^{(q)}_{i})=n-p^q.$ Moreover, we have

\begin{proposition}\label{p2.1} Let  $|\omega|=n$ be even and $i\in Y_0 $. Then
 $\Phi^{(q)}_{i}\in {\rm Der}(\mathcal{H},\mathcal{W}).$
\end{proposition}
  \begin{proof} It is sufficient to verify the following
equation for $\mathrm{D_{H}}(x^{(\alpha)}x^u),$
$\mathrm{D_{H}}(x^{(\beta)}x^v)\in \mathcal{H}:$
\begin{eqnarray}
&&\Phi^{(q)}_{i}([\mathrm{D_{H}}(x^{(\alpha)}x^u),\mathrm{D_{H}}(x^{(\beta)}x^v)])
\nonumber\\
&=& [\Phi^{(q)}_{i}
(\mathrm{D_{H}}(x^{(\alpha)}x^u)),\mathrm{D_{H}}(x^{(\beta)}x^v)]
+[\mathrm{D_{H}}(x^{(\alpha)}x^v),
\Phi^{(q)}_{i}(\mathrm{D_{H}}(x^{(\beta)}x^u))].\label{e6.6.20}
\end{eqnarray}
The verification is divided into three parts.\newline

 \noindent\textit{Case} (i): $u=v=\emptyset.$ The left-hand side of
(\ref{e6.6.20}) is as follows:
\begin{equation}
 \Phi_{i}^{(q)}([\mathrm{D_{H}}(x^{(\alpha)}),\ \mathrm{D_{H}}(x^{(\beta)})])
=\Phi_{i}^{(q)}(\mathrm{D_{H}}(\mathrm{D_{H}}(x^{(\alpha)})(x^{(\beta)}))=
\partial_{i}^{p^{q}}(\mathrm{D_{H}}(x^{(\alpha)})(x^{(\beta)}))
\mathrm{D_{H}}(x^{\omega}).\label{e6.6.21}
\end{equation}
By (\ref{le0}), the right-hand side equals:
\begin{eqnarray}
&&[x^{(\alpha-p^{q}\varepsilon_i)}\mathrm{D_{H}}(x^{\omega})
,\mathrm{D_{H}}(x^{(\beta)})] + [\mathrm{D_{H}}(x^{(\alpha)}), \
x^{(\beta-p^{q}\varepsilon_i)}
\mathrm{D_{H}}(x^{\omega})]\nonumber
\\
 &=& \big(-\mathrm{D_{H}}(x^{(\beta)})(x^{(\alpha-p^{q}\varepsilon_i)})
  +\mathrm{D_{H}}(x^{(\alpha)})(x^{(\beta-p^{q}\varepsilon_i)})\big)
  \mathrm{D_{H}}(x^{\omega}).\label{e6.6.22}
 \end{eqnarray}
 Noticing that $\partial_{i}^{p^{q}}$ is a derivation, one can compute
 the coefficient of $\mathrm{D_{H}}(x^{\omega})$ in (\ref{e6.6.21}):
\begin{eqnarray*}
\partial_i^{p^{q}}\big(\mathrm{D_{H}}(x^{(\alpha)})(x^{(\beta)})\big)&=&
  \partial_i^{p^{q}} \Big(\sum_{r\in Y_0}\tau(r)
  x^{(\alpha-\varepsilon_r)}\partial_{r'}(x^{(\beta)})\Big)\nonumber\\
  &=&\sum_{r\in Y_0}\big(\tau(r)
  x^{(\alpha-\varepsilon_r-p^q\varepsilon_i)}x^{(\beta-\varepsilon_{r'})}
  +\tau(r)
  x^{(\alpha-\varepsilon_r)}x^{(\beta-\varepsilon_{r'}-p^{q}\varepsilon_i)}\big)
 \end{eqnarray*}
 The coefficient of $\mathrm{D_{H}}(x^{\omega})$ in (\ref{e6.6.22}) is as
 follows:
\begin{eqnarray*}
-\mathrm{D_{H}}(x^{(\beta)})(x^{(\alpha-p^{q}\varepsilon_i)})
&+&\mathrm{D_{H}}(x^{(\alpha)})(x^{(\beta-p^{q}\varepsilon_i)})\\
&=&-\sum_{r\in Y_0}\tau(r')
  x^{(\beta-\varepsilon_{r'})}x^{(\alpha-\varepsilon_{r}-p^{q}\varepsilon_{i})}
  +\sum_{r\in Y_0}\tau(r)
  x^{(\alpha-\varepsilon_r)}x^{(\beta-\varepsilon_{r'}-p^{q}\varepsilon_i)}\\
  &=&\sum_{r\in Y_0}(\tau(r)x^{(\alpha-\varepsilon_{r}-p^{q}\varepsilon_{i})}
   x^{(\beta-\varepsilon_r')}+\tau(r)x^{(\alpha-\varepsilon_r)}
     x^{(\beta-\varepsilon_{r'}-p^{q}\varepsilon_{i})}).
 \end{eqnarray*}%(5)
 Hence (\ref{e6.6.20}) holds in this case.\newline

 \noindent\textit{Case} (ii): $u\not= \emptyset$ and
  $v\not= \emptyset.$ Then $|u|\geq 2,$
  $|v|\geq 2.$ By the definition of $\Phi_{i}^{(q)},$ it is easily
  seen that the two sides of (\ref{e6.6.20}) are all zero.\newline

 \noindent\textit{Case} (iii): $u=\emptyset$
  $v\not= \emptyset. $
By the definition of $\Phi_{i}^{(q)},$ the left-hand side is as
follows:
\begin{eqnarray*}
 \Phi_{i}^{(q)}([\mathrm{D_{H}}(x^{(\alpha)}),
 \mathrm{D_{H}}(x^{(\beta)}x^v)])
 &=& \Phi_{i}^{(q)}(\mathrm{D_{H}}(\mathrm{D_{H}}(x^{(\alpha)})(x^{(\beta)}x^v)))
 \\
 &=&
 \partial_{i}^{p^{q}}(\mathrm{D_{H}}(x^{(\alpha)})(x^{(\beta)}x^v))\mathrm{D_{H}}(x^{\omega})
 = 0.\quad (\mbox{Note}\ |v|\geq 2)
 \end{eqnarray*}
The right-hand side is as follows:
\begin{eqnarray*}
 & &[\Phi_{i}^{(q)}(\mathrm{D_{H}}(x^{(\alpha)})),\ \mathrm{D_{H}}(x^{(\beta)}x^v)]
+
 [\mathrm{D_{H}}(x^{(\alpha)}),\Phi_{i}^{(q)}(\mathrm{D_{H}}(x^{(\beta)}x^v))]
\\
&=&[\partial^{p^{q}}_i(x^{(\alpha)}) \mathrm{D_{H}}(x^{\omega}), \
\mathrm{D_{H}}(x^{(\beta)}x^v)] +
 [\mathrm{D_{H}}(x^{(\alpha)}),\ \partial^{p^{q}}_i(x^{(\beta)}x^v) (\mathrm{D_{H}}(x^{\omega}))]\\
 &=&[\partial^{p^{q}}_i(x^{(\alpha)}) \mathrm{D_{H}}(x^{\omega}),\
 \mathrm{D_{H}}(x^{(\beta)}x^v)]\\
 &=&-\mathrm{D_{H}}(x^{(\beta)}x^v)(\partial^{p^{q}}_i(x^{(\alpha)})\mathrm{D_{H}}(x^{\omega}))+
\partial^{p^{q}}_i(x^{(\alpha)})[\mathrm{D_{H}}(x^{\omega}), \mathrm{D_{H}}(x^{(\beta)}x^v)]\\
&=&-\mathrm{D_{H}}(x^{(\beta)}x^v)(\partial^{p^{q}}_i(x^{(\alpha)}))\mathrm{D_{H}}(x^{\omega})\\
&=& 0.
 \end{eqnarray*}
The proof is complete.
\end{proof}

Now we define the second series of exceptional
derivations.   We have known that $({\rm ad}\partial_i)^{p^{q}}$
is a derivation of $\mathcal{H}.$ Define for $i\in Y_{0}$ and
$q\in \mathbb{N},$
$$\Theta^{(q)}_{i}: \mathcal{H} \rightarrow \mathcal{W},\quad
 \mathrm{D_{H}}(f)\mapsto  x^{\omega}\,
 (\mathrm{ad}\partial_i)^{p^q}(\mathrm{D_{H}}(f)).
$$ By the definition, $\mathrm{zd}(\Theta^{(q)}_{i})=n-p^q.$
Clearly, $\Theta^{(q)}_{i}$ may be naturally extended to  a
linear mapping of $\mathcal{W}.$ Note that
$\Theta^{(q)}_{i}(\mathrm{D_{H}}(f))=x^{\omega}\mathrm{D_{H}}(\partial_{i}^{p^{q}}(f)).$ Moreover, we have the following

\begin{proposition} \label{p2.2}
Let $|\omega|=n$ be even and $i\in Y_0$. Then
 $ \Theta^{(q)}_{i}\in {\rm Der}(\mathcal{H},\mathcal{W}).$
\end{proposition}
 \begin{proof} View $\Theta^{(q)}_{i}$
 as the
linear mapping of $\mathcal{W} $  and consider the standard basis elements
$x^{(\alpha)}x^u\partial_r$ of $\mathcal{W}$
 satisfying that $r\not\in u $ and that if $u=\emptyset$ then $r\in Y_{0}.$
   Since every element of
$\mathrm{D_{H}}(f) \in\mathcal{H}$ must be a linear combination of
such standard basis elements, it suffices to show that for such
standard basis elements $x^{(\alpha)}x^u\partial_r$ and
$x^{(\beta)}x^v\partial_s,$ the following holds:
\begin{eqnarray}
&&\Theta^{(q)}_{i}([x^{(\alpha)}x^u\partial_r,\
x^{(\beta)}x^v\partial_s] )\nonumber\\
&=&[\Theta^{(q)}_{i}(x^{(\alpha)}x^u\partial_r),x^{(\beta)}x^v\partial_s] +
 [x^{(\alpha)}x^u\partial_r,\ \Theta^{(q)}_{i}(x^{(\beta)}x^v\partial_s)
 ].\label{e6.6.23}
\end{eqnarray}
When $u=v=\phi,$ it is easy to see that (\ref{e6.6.23}) holds.
When
 $u\not=\emptyset$ and
 $v\not=\emptyset,$ the two sides of (\ref{e6.6.23})  vanish. It remains only the case that
 $u=\emptyset$ but
 $v\not=\emptyset.$ Note that $r\in Y_{0}$ in this case.
 Consequently, the left-hand side of (\ref{e6.6.23}) vanishes. As
 $s\not\in v,$ the first summand in the right-hand side  of (\ref{e6.6.23}) is zero.
 Clearly, the second is also zero.
 \end{proof}

\begin{remark}\label{r2.3}
From the proof of Proposition \ref{p2.2}, it is easily seen that
in general, $\Theta^{(q)}_{i} \not\in {\rm Der}\mathcal{W}$ and
 $\Theta^{(q)}_{i}\not\in {\rm Der}(\mathcal{S},\mathcal{W}) $
  where $\mathcal{S}$ denotes the even part of
  the special superalgebra and  $\Theta^{(q)}_{i}$ is naturally extended.
\end{remark}

 Let us consider the third series of exceptional derivations.
 Define for $i\in Y_0,$
 $$\Psi^{(i)}: \mathcal{H}\longrightarrow \mathcal{W},\quad
\mathrm{D_{H}}(f)\longmapsto
 \partial_i\partial_{i'}(f)\mathrm{D_{H}}(x
^{\omega})\quad \mbox{for}\,\, f\in \mathcal{O}(m,n;\underline{t}).
$$
As $\ker(\mathrm{D_{H}})=\mathbb{F}\, 1,$
 the linear mapping $\Psi^{(i)}$ is well defined and
 ${\rm zd}(\Psi^{(i)})=n-2.$ Moreover, we have the
 following

\begin{proposition}\label{p2.4}
 Let  $|\omega|=n$ be even and $i\in Y_0 $. Then
 $\Psi^{(i)}\in {\rm Der}(\mathcal{H},\mathcal{W}).$
\end{proposition}
 \begin{proof}
 We want to verify that for $\alpha, \beta \in \mathbb{A},$ $u, v\in \mathbb{B}^0,$
 the following holds:
\begin{eqnarray}
&&\Psi^{(i)}([\mathrm{D_{H}}(x^{(\alpha)}x^u),\mathrm{D_{H}}(x^{(\beta)}x^v)])\nonumber\\
&=&
[\Psi^{(i)}(\mathrm{D_{H}}(x^{(\alpha)}x^u)),\mathrm{D_{H}}(x^{(\beta)}x^v)]
+[\mathrm{D_{H}}(x^{(\alpha)}x^u),
\Psi^{(i)}(\mathrm{D_{H}}(x^{(\beta)}x^v))].\label{e6.6.24}
\end{eqnarray}

 \noindent\textit{Case} (i): $u\not=\emptyset,$ $v\not=\emptyset.$
Since $|u|\geq 2,$ $|v|\geq 2 $ in this case, two sides of (\ref{e6.6.24})
vanish.\newline

  \noindent\textit{Case} (ii):
$u\not=\emptyset,$ $v=\emptyset.$ Then the right-hand side of
(\ref{e6.6.24}) is as follows:
\begin{eqnarray*}
&&[\mathrm{D_{H}}(x^{(\alpha)}x^u),\
x^{(\beta-\varepsilon_i-\varepsilon_{i'})}
\mathrm{D_{H}}(x^{\omega})]\\
&=&\mathrm{D_{H}}(x^{(\alpha)}x^u)(x^{(\beta-\varepsilon_i-\varepsilon_{i'})})
\mathrm{D_{H}}(x^{\omega})
+x^{(\beta-\varepsilon_i-\varepsilon_{i'})}
[\mathrm{D_{H}}(x^{(\alpha)}x^u),\ \mathrm{D_{H}}(x^{\omega})]\\
&=&0.
\end{eqnarray*}
The left-hand side is as follows
\begin{eqnarray*}
 \Psi^{(i)}(\mathrm{D_{H}}(\mathrm{D_{H}}(x^{(\alpha)}x^u)(x^{(\beta)}))=
\Psi^{(i)}\Big(\mathrm{D_{H}}\Big(\sum_{r\in
Y_0}\tau(r)x^{(\alpha-\varepsilon_r)}x^ux^{(\beta-r')}\Big)\Big)=0.
\end{eqnarray*}
This proves (\ref{e6.6.24}) in this case.\newline

 \noindent\textit{Case} (iii): $u=v=\phi.$ The right-hand side of (\ref{e6.6.24}) is
as follows:
\begin{eqnarray*}
&&[x^{(\alpha-\varepsilon_i-\varepsilon_{i'})}\mathrm{D_{H}}(x^{\omega}),\
\mathrm{D_{H}}(x^{(\beta)})]
 +[\mathrm{D_{H}}(x^{(\alpha)}), \
 x^{(\beta-\varepsilon_i-\varepsilon_{i'})}\mathrm{D_{H}}(x^{\omega})]\\
 &=& -\mathrm{D_{H}}(x^{(\beta)})(x^{(\alpha-\varepsilon_i-\varepsilon_{i'})})
 \mathrm{D_{H}}(x^{\omega})
 +\mathrm{D_{H}}(x^{(\alpha)})(x^{(\alpha-\varepsilon_i-\varepsilon_{i'})})
 \mathrm{D_{H}}(x^{\omega}).
\end{eqnarray*}
The left-hand side is as follows:
\begin{eqnarray*}
&& \Psi^{(i)}(\mathrm{D_{H}}(\mathrm{D_{H}}(x^{(\alpha)})(x^{(\beta)}))\\
&=&\partial_i\partial_{i'}(\mathrm{D_{H}}(x^{(\alpha)})(x^{(\beta)}))\mathrm{D_{H}}(x^{\omega})\\
&=&\partial_i\partial_{i'}\Big(\sum_{r\in
Y_0}\tau(r)(x^{(\alpha-\varepsilon_{r})}x^{(\beta-
\varepsilon_{r'})}\Big)\mathrm{D_{H}}(x^{\omega})\\
&=&\partial_i\Big(\sum_{r\in
Y_0}\tau(r)\big(x^{(\alpha-\varepsilon_{r}-\varepsilon_{i'})}x^{(\beta-
\varepsilon_{r'})}+x^{(\alpha-\varepsilon_{r})}
x^{(\beta-\varepsilon_{r'}-\varepsilon_{i'})}\big)\Big)\mathrm{D_{H}}(x^{\omega})\\
&=&\sum_{r\in
Y_0}\tau(r)\big(x^{(\alpha-\varepsilon_{r}-\varepsilon_{i'}-\varepsilon_{i})}x^{(\beta-
\varepsilon_{r'})}+ x^{(\alpha-\varepsilon_{r}-\varepsilon_{i'})}
x^{(\beta-\varepsilon_{r'}-\varepsilon_{i})} \\&&+
x^{(\alpha-\varepsilon_{r}-\varepsilon_{i})}
x^{(\beta-\varepsilon_{r'}-\varepsilon_{i'})}+
x^{(\alpha-\varepsilon_{r})}
x^{(\beta-\varepsilon_{r'}-\varepsilon_{i'}-\varepsilon_{i})}
 \big)\mathrm{D_{H}}(x^{\omega})\\
 &=&\sum_{r\in
Y_0}\tau(r)x^{(\alpha-\varepsilon_{r}-\varepsilon_{i'}-\varepsilon_{i})}x^{(\beta-
\varepsilon_{r'})}\mathrm{D_{H}}(x^{\omega})+
  \Big(\sum_{r\in Y_0}\tau(r)(x^{(\alpha-\varepsilon_{r})}
x^{(\beta-\varepsilon_{r'}-\varepsilon_{i'}-\varepsilon_{i})})\Big)
\mathrm{D_{H}}(x^{\omega})\\
&=& -\Big(\sum_{r\in
Y_0}\tau(r')x^{(\alpha-\varepsilon_{r}-\varepsilon_{i'}-\varepsilon_{i})}x^{(\beta-
\varepsilon_{r'})} \Big)\mathrm{D_{H}}(x^{\omega})+
  \Big(\sum_{r\in Y_0}\tau(r)\big(x^{(\alpha-\varepsilon_{r})}
x^{(\beta-\varepsilon_{r'}-\varepsilon_{i'}-\varepsilon_{i})}\big)\Big)
\mathrm{D_{H}}(x^{\omega})
\\
&=&\big(-\mathrm{D_{H}}\big(x^{(\beta)}\big)
\big(x^{(\alpha-\varepsilon_{i'}-\varepsilon_{i})}\big)
+\mathrm{D_{H}}\big(x^{(\alpha)}\big)\big(x^{(\beta-\varepsilon_{i'}
-\varepsilon_{i})}\big)\big)\mathrm{D_{H}}(x^{\omega}).
\end{eqnarray*}
Summarizing, (\ref{e6.6.24}) holds, completing the proof.
\end{proof}

 \section{$\mathbb{Z}$-homogeneous derivations}

 In this section, we first determine the derivations from $\frak{N}$ into $\mathcal{W}$ which vanish
on the top of $\frak{N}$. To that aim, one needs to investigate the action on the
generators of $\frak{N}$ for such a derivation. Recall the generating  theorem established in Section 2. We shall consider
the set $\mathcal{M}$ and $\mathcal{N}$ separately. Recall
$$\mathcal{G}=C_{\mathcal{W} } (\mathcal{W} _{-1})
=\mathrm{span}_\mathbb{F}\{x^u\partial_r\,|\,r\in Y,u\in
\mathbb{B},\mathrm{p}(x^u\partial_r)=\overline{0}\}.$$ For simplicity,
 put $$E(\mathcal{G}):=\oplus _{r\in \mathbb{Z}}\mathcal{G}_{2r}, \quad
O(\mathcal{G}):=\oplus _{r\in \mathbb{Z}}\mathcal{G}_{2r+1}.$$
We shall frequently use the following simple fact.
 \begin{lemma}\label{bc2.2} Suppose $\phi\in {\rm Der}(\frak{N},\mathcal{W})$ satisfies $\phi(\frak{N}_{[-1]})=0.$
Then for $D\in \frak{N},$
$\phi([D,\frak{N}_{[-1]}])=0$ if and only if $\phi(D)\in \mathcal{G}.$
\end{lemma}

Let us first consider the elements in $\mathcal{M}.$
\begin{lemma}\label{L3.1}
Let $\phi\in {\rm Der}(\frak{N},\mathcal{W})$ be homogeneous such
that $\phi(\frak{N}_{[-1]}\oplus \frak{N}_{[0]})=0.$ Suppose
${\rm zd}(\phi)+a$ is odd and
$\phi(\mathrm{D_{H}}(x^{(b\varepsilon_i)}))=0$ for all $b<a,$
where $a\leq\pi_i$ is a fixed positive integer  and $i\in Y_0.$
Then the following statements hold.

{\rm{(i)}}\;\, If $a\not\equiv 1\pmod{p}$ then
$\phi(\mathrm{D_{H}}(x^{(a\varepsilon_i)}))=0.$

{\rm{(ii)}}\;\ If $a\equiv 1\pmod{p}$ and $a-1$ is not any
$p$-power, then $\phi(\mathrm{D_{H}}(x^{(a\varepsilon_i)}))=0.$

{\rm{(iii)}} If $a-1=p^q$ for some $q\in \mathbb{N},$ then there are
$\mu^{(q)}_i, \eta^{(q)}_i\in \mathbb{F}$ such that
$\mu^{(q)}_i \eta^{(q)}_i=0$ and
$$(\phi-\tau(i)\mu^{(q)}_i\Theta^{(q)}_{i}-
\tau(i)\eta^{(q)}_i ({\rm
ad}\partial_i)^{p^q})(\mathrm{D_{H}}(x^{(a\varepsilon_i)}))=0.$$
\end{lemma}
 \begin{proof}  Since $\phi(\frak{N}_{[-1]}) =0$ and
${\rm zd}(\phi)+a$ is odd, by Lemma \ref{bc2.2},
$\phi(\mathrm{D_{H}}(x^{(a\varepsilon_i)}))\in O(\mathcal{G}).$
Thus one may assume that
$$\phi(\mathrm{D_{H}}(x^{(a\varepsilon_i)}))=\sum_{r\in Y_0} f_r \partial_r \quad
{\rm where }\  f_r \in \Lambda (n).$$ For arbitrary $j\in
Y_0\setminus \{i,i'\},$ we have $[\mathrm{D_{H}}(x_jx_{j'}), \
\mathrm{D_{H}}(x^{(a\varepsilon_i)})]=0.$ Applying $\phi $ to this
equation, one gets
$$[\tau(j)x_{j'}\partial_{j'}+\tau(j')x_j\partial_j,\ \sum_{r\in Y_0} f_r \partial_r]=0.$$
It follows that $f_j=0$  for all $j\in Y_0\setminus \{i,i'\}$ and
therefore,
$$\phi(\mathrm{D_{H}}(x^{(a\varepsilon_i)}))=f_i\partial_i +f_{i'}\partial_{i'}.$$
Similarly, applying $\phi$ to the equation
$[\mathrm{D_{H}}(x^{(2\varepsilon_i)}),\mathrm{D_{H}}(x^{(a\varepsilon_i)})]=0,$
one gets $f_i=0.$ Thus
$$\phi(\mathrm{D_{H}}(x^{(a\varepsilon_i)}))=f_{i'}\partial_{i'}.$$
For arbitrary $k,l\in Y_1$ with $k\not=l,$ we have
$$0=[\mathrm{D_{H}}(x_kx_l),\phi(\mathrm{D_{H}}(x^{(a\varepsilon_i)}))]
=(x_l\partial_k-x_k\partial_l)(f_{i'})\partial_{i'}.$$
Hence $x_l\partial_k(f_{i'})-x_k\partial_l(f_{i'})=0$ and therefore,
$x_k\partial_l(f_{i'})=0.$ This implies that $f_{i'}\in\mathbb{ F}\,
x^{\omega} $ or $f_{i'}\in\mathbb{ F}.$ Therefore, there is
$\mu_i\in\mathbb{ F},$ such that

\begin{equation}\phi(\mathrm{D_{H}}(x^{(a\varepsilon_i)}))=\mu_ix^{\omega}\partial_{i'}\quad
{\rm or} \quad \mu_i\partial_{i'}.\label{e6.6.25}
\end{equation}%(3.1)
Suppose  $a-1=p^{q}$ for some $q\in \mathbb{N}.$ Then it is
easily seen that (iii) holds. Note that
$[\mathrm{D_{H}}(x_ix_{i'}),\mathrm{D_{H}}(x^{(a\varepsilon_i)})]
=\tau(i')a\mathrm{D_{H}}(x^{(a\varepsilon_i)}).$ Applying $\phi$,
we have
$$[\tau(i)x_{i'}\partial_{i'}+\tau(i')x_i\partial_i,\ \phi(\mathrm{D_{H}}(x^{(a\varepsilon_i)}))]=
\tau(i')a\phi(\mathrm{D_{H}}(x^{(a\varepsilon_i)})).$$ It follows
from (\ref{e6.6.25}) that $(a-1)\mu_i=0.$ Consequently, (i) holds.

To prove (ii),  suppose  $a\equiv 1\pmod p$ and $a-1$ is not
any $p$-power. Write $a-1$ to be the $p$-adic  form
$a-1=\sum^{t}_{r=1}c_rp^r$ where $c_t\not=0.$ Then
$\binom{a}{p^t}\not\equiv 0\pmod p.$ Note that
\begin{equation}[\mathrm{D_{H}}(x^{(p^t\varepsilon_i)}x_{i'}),
\mathrm{D_{H}}(x^{((a-p^t+1)\varepsilon_i)})]=
\tau(i')\binom{a}{p^t}\mathrm{D_{H}}(x^{(a\varepsilon_i)}).
\label{e6.6.26}
\end{equation}
It is clear that $a-p^t+1< a.$ Then
$\phi(\mathrm{D_{H}}(x^{((a-p^t+1)\varepsilon_i)}))=0.$ We want to
show that $\phi(\mathrm{D_{H}}(x^{(p^t\varepsilon_i)}x_{i'}))=0.$
Since $p^t <a-1,$ it follows that
$$\phi(\mathrm{D_{H}}(x^{(p^{t}\varepsilon_i)}x_{i'}))
=\tau(i')[\mathrm{D_{H}}(x^{(2\varepsilon_{i'})}),\
\phi(\mathrm{D_{H}}(x^{((p^{t}+1)\varepsilon_{i})}))]=0.$$ Thus  we obtain
from (\ref{e6.6.26}) that
$\phi(\mathrm{D_{H}}(x^{(a\varepsilon_{i})}))=0;$ that is, (ii)
holds.
\end{proof}

\begin{lemma}\label{L3.2}
 Let $\phi\in {\rm Der}(\frak{N},\mathcal{W})$ be homogeneous such
that $\phi(\frak{N}_{[-1]}\oplus \frak{N}_{[0]})=0.$ Suppose
$a\in \mathbb{N}$ and $i\in Y_0$ such that
$\phi(\mathrm{D_{H}}(x^{(b\varepsilon_i)}))=0$ for all $b<a.$ If
${\rm zd}(\phi)+a$ is even, then the following statements holds.

{\rm{(i)}}\;\, If $a\not\equiv 0\pmod{p}$ then
$\phi(\mathrm{D_{H}}(x^{(a\varepsilon_i)}))=0.$

{\rm{(ii)}}\; If $a\equiv 0\pmod{p}$ and $p$ is not any $p$-power,
then $\phi(\mathrm{D_{H}}(x^{(a\varepsilon_i)}))=0.$

{\rm{(iii)}} If $a=p^q$ for some $q\in \mathbb{N},$ then there is
$\lambda^{(q)}_{i}\in\mathbb{F}$ such that

$$(\phi-\lambda^{(q)}_{i}\Phi^{(q)}_{i})(\mathrm{D_{H}}(x^{(a\varepsilon_i)}))=0.
$$
\end{lemma}
 \begin{proof} Since $\phi(\frak{N}_{[-1]}) =0$ and
${\rm zd}(\phi)+a$ is even, by Lemma \ref{bc2.2},
$\phi(\mathrm{D_{H}}(x^{(a\varepsilon_i)}))\in E(\mathcal{G}).$
Thus we may assume that
\begin{equation}\phi(\mathrm{D_{H}}(x^{(a\varepsilon_i)}))
=\sum_{r\in Y_1} f_r \partial_r \quad {\rm where }\  f_r \in \Lambda
(n).\label{e6.6.27}
\end{equation}

(i) Applying $\phi$ to the equation
$$[\mathrm{D_{H}}(x_ix_{i'}),\ \mathrm{D_{H}}(x^{(a\varepsilon_i)})]
=\tau(i')a\mathrm{D_{H}}(x^{(a\varepsilon_i)})$$ one may obtain
$\phi(\mathrm{D_{H}}(x^{(a\varepsilon_i)}))=0,$ since $a\not\equiv
0\pmod p.$

(ii) Write $a$ to be the $p$-adic form $a=\sum^t_{r=1}c_rp^r,\
c_t\neq 0.$ Then $\binom{a}{p^t}\not\equiv 0\pmod p.$ Just as in
the proof of Lemma \ref{L3.1} we have
\begin{equation}
[\mathrm{D_{H}}(x^{(p^t\varepsilon_i)}x_{i'}),
\mathrm{D_{H}}(x^{((a-p^t+1)\varepsilon_i)})]=
\tau(i')\binom{a}{p^t}\mathrm{D_{H}}(x^{(a\varepsilon_i)}).\label{e6.6.28}
\end{equation}%3.4
Clearly, $\phi(\mathrm{D_{H}}(x^{((a-p^t+1)\varepsilon_i)}))=0,$
since $a-p^t+1<a.$ On the other hand, we  have also
\begin{equation*}
[\mathrm{D_{H}}(x^{(2\varepsilon_{i'})}),
\mathrm{D_{H}}(x^{(a\varepsilon_i)})]=
\tau(i')\mathrm{D_{H}}(x^{(a-1)\varepsilon_i}x_{i'}).
\end{equation*}
It follows from (\ref{e6.6.27}) and the equation above that
$\phi(\mathrm{D_{H}}(x^{(a-1)\varepsilon_i}x_{i'}))=0.$ Then
$$
\phi(\mathrm{D_{H}}(x^{(p^t\varepsilon_i)}x_{i'}))=0,
$$
since $p^t<a-1.$ Now (ii) follows from (\ref{e6.6.28}).

(iii) For $k,l\in Y_1$ with $k\not=l,$ we have
\begin{equation}
0=[\mathrm{D_{H}}(x_kx_l),\phi(\mathrm{D_{H}}(x^{(a\varepsilon_i)}))]=\sum_{r\in
Y_1}(x_l\partial_k-x_k\partial_l)(f_{r})\partial_{r}-f_l\partial_k+f_k\partial_l.\label{e6.6.29}%3.5
\end{equation}
It follows that $(x_l\partial_k-x_k\partial_l)(f_{r})=0$ for $k,l\in
Y_1\setminus \{r\}.$  This implies that $x_k\partial_l(f_{r})=0$ whenever
$k,l\in Y_1\setminus\{r\},\ k\not=l.$ From this one may deduce
that
$$f_r\in \Lambda(n)_{[n-1]}\cup \Lambda(n)_{[n]}\cup
\mathbb{F}\quad\mbox{for all}\ r\in Y_1.
$$
Clearly, $f_r$ cannot be nonzero in $\mathbb{F},$ since ${\rm
zd}(\phi) +a$  is even. On the other hand, if $f_r\in
\Lambda(n)_{[n]},$ one may obtain $f_r=0$ from (\ref{e6.6.29}).
Summarizing, $f_r\in\Lambda(n)_{[n-1]}.$ This implies that if $n$
is odd then $f_r=0 $ for all  $r\in Y_1.$ Thus we assume that $n$
is even in the following. From the fact that $x_k\partial_l(f_r)=0 $
whenever $k,l\in Y_1\setminus\{r\},\ k\not=l,$ one may deduce that
$f_r\in \mathbb{F}\, x^{\omega-<r>}.$ Hence we may assume that
\begin{equation}
\phi(\mathrm{D_{H}}(x^{(a\varepsilon_i)}))=\sum_{r\in
Y_1}c_rx^{\omega-<r>}\partial_r\quad {\rm where} \quad c_r\in \mathbb{F}.
\label{e6.6.30}%3.6
\end{equation}
From (\ref{e6.6.29}) we have
$$-x_k\partial_l(c_kx^{\omega-<k>})\partial_k-c_lx^{\omega-<l>}\partial_k+
x_l\partial_k(c_lx^{\omega-<l>})\partial_l+c_kx^{\omega-<k>}\partial_l=0.$$
Furthermore,
$$x_k\partial_l(c_kx^{\omega-<k>})+c_lx^{\omega-<l>}=0.$$
Without loss of generality, we may assume that $k<l.$ Then the
equation
yields that $(-1)^{l-2}(-1)^{k-1}c_k+c_l=0;$ that is,
$c_k=(-1)^{k+l}c_l.$ Recall that $a=p^q.$ Put
$\lambda^{(q)}_{i}:=c_{2m+1}.$ Then we obtain from (\ref{e6.6.30})
that
$$\phi(\mathrm{D_{H}}(x^{(a\varepsilon_i)}))
=\lambda^{(q)}_{i}\mathrm{D_{H}}(x^{\omega}).$$
Therefore,
$$(\phi- \lambda^{(q)}_{i}\Phi^{(q)}_{i})(\mathrm{D_{H}}(x^{(a\varepsilon_i)}))=0.$$
The proof is complete.
\end{proof}

Now we can reduce the derivations vanishing on the top of $\frak{N}$ to be vanishing on $\mathcal{M}.$
\begin{lemma}\label{L3.4}
Let $\phi\in {\rm Der}(\frak{N},\mathcal{W}).$ Suppose   ${\rm
zd}(\phi) $ is even and $\phi(\frak{N}_{[-1]}\oplus
\frak{N}_{[0]})=0.$ Then
$\phi(\mathrm{D_{H}}(x^{(a_i\varepsilon_i)}))=0$ for all $i\in Y_0
$ and $a_i\leq \pi_i.$ That is, $\phi(\mathcal{M})=0.$
\end{lemma}
 \begin{proof} We proceed by induction on $a_i$. Assume
that the assertion holds for $a_i-1.$ Let us consider the case
$a_i$. If $a_{i}=p^q$ for some $q\in \mathbb{N},$ then ${\rm
zd}(\phi)+a_{i} $ is odd. By Lemma \ref{L3.1}(i),
$\phi(\mathrm{D_{H}}(x^{(a_i\varepsilon_i)}))=0.$ Assume that
$a_i$ is not any $p$-power. The discuss is divided into the
following two cases (i) $a_i\equiv 0\pmod p$ and (ii)
$a_i\not\equiv 0\pmod p$. In the case (i), if ${\rm
zd}(\phi)+a_{i} $ is odd, then Lemma \ref{L3.1}(i) ensures that
$\phi(\mathrm{D_{H}}(x^{(a_i\varepsilon_i)}))=0;$ if ${\rm
zd}(\phi)+a_{i} $ is even, then Lemma \ref{L3.2}(ii) ensures that
$\phi(\mathrm{D_{H}}(x^{(a_i\varepsilon_i)}))=0,$ since $a_{i}$ is
not any $p$-power. In the case (ii), if ${\rm zd}(\phi)+a_{i} $ is
even, then Lemma \ref{L3.2}(i) ensures that
$\phi(\mathrm{D_{H}}(x^{(a_i\varepsilon_i)}))=0;$ if ${\rm
zd}(\phi)+a_{i} $ is odd, then $a_{i}$ is odd and there is no $q$
such that $a_{i}-1=p^q $ and therefore,
$\phi(\mathrm{D_{H}}(x^{(a_i\varepsilon_i)}))=0$ by Lemma \ref{L3.1}.
  The proof is complete.
  \end{proof}

\begin{lemma}\label{L3.5}
Let $\phi\in {\rm Der}(\frak{N},\mathcal{W})$ and  ${\rm zd}(\phi)
$ be odd. Suppose
 $\phi(\frak{N}_{[-1]}\oplus \frak{N}_{[0]})=0.$ Then there are
$\lambda^{(s_r)}_r, \mu^{(s_r)}_r, \eta^{(s_r)}_r\in \mathbb{F},$
where $r\in Y_0$ and $1\leq s_r < t_r, $ such that
$\mu^{(s_r)}_r\eta^{(s_r)}_r=0$ and
$$\Big(\phi-\sum_{r\in Y_0}\sum_{s_r=1}^{t_r-1}\big(\lambda^{(s_r)}_r\Phi^{(s_r)}_r
+\mu^{(s_r)}_r\Theta^{(s_{r})}_{r}+\eta^{(s_r)}_r({\rm
ad}\partial_r)^{p^{s_r}}\big)\Big)(\mathrm{D_{H}}(x^{(a_i\varepsilon_i)}))=0$$
for all $i\in Y_0,$ $1\leq a_i\leq \pi_i.$
\end{lemma}
 \begin{proof} For $1<a_i<p,$ by virtue of Lemmas
\ref{L3.1}(i)
 and \ref{L3.2}(i),
 one may show by
induction on $a_i$ that
$$\phi(\mathrm{D_{H}}(x^{a_i\varepsilon_i)}))=0,\quad 1<a_{i}<p,\ i\in Y_0.$$
Put $\phi^{(0)}:=\phi.$ Assume  inductively that we have
constructed a derivation $\phi^{(q)}$ from $\frak{N}$ into $\mathcal{W}$ for $q\geq 0$ such that
$$\phi^{(q)}(\mathrm{D_{H}}(x^{a_i\varepsilon_i)}))=0 \quad \mbox{ for all}
\  i\in Y_0 \ {\rm and}\  1<a_{i}< p^{q+1}.$$
 We first consider $\phi^{(q)}(\mathrm{D_{H}}(x^{(p^{q+1}\varepsilon_i)})).$
 Note that ${\rm zd}(\phi)+ p^{q+1}$ is even. By Lemma \ref{L3.2}(iii),
  there is $\lambda^{(q+1)}_r\in \mathbb{F}$ for $r\in Y_0$ such that
 $$(\phi^{(q)}-\lambda^{(q+1)}_r\Phi^{(q+1)}_r)
 (\mathrm{D_{H}}(x^{(p^{q+1}\varepsilon_r)}))=0\quad
 {\rm for} \  r\in Y_0.$$
 Put $\widetilde{\phi}^{(q+1)}:=
\phi^{(q)}-\sum_{r\in Y_0 }\lambda^{(q+1)}_r\Phi^{(q+1)}_r.$ Then
$$\widetilde{\phi}^{(q+1)}(\mathrm{D_{H}}(x^{(p^{q+1}\varepsilon_r)}))=0
\quad {\rm for\  all} \  r\in Y_0.$$ By Lemma \ref{L3.1}(iii),
 there are $\mu^{(q+1)}_r, \eta^{(q+1)}_r\in \mathbb{F}$ such that
$\mu^{(q+1)}_r\eta^{(q+1)}_r=0$ and
$$(\widetilde{\phi}^{(q+1)}-\tau(r)\mu^{(q+1)}_r\Theta_{i}^{(q+1)}- \tau(r)\eta^{(q+1)}_r({\rm
ad}\partial_r)^{p^{q+1}}
)(\mathrm{D_{H}}(x^{((p^{q+1}+1)\varepsilon_r)}))=0,\ r\in Y_0.$$
Put $\phi^{(q+1)}:=\widetilde{\phi}^{(q+1)}-\sum_{r\in
Y_0}\tau(r)\mu^{(q+1)}_r\Theta_{r}^{(q+1)} -\sum_{r\in
Y_0}\tau(r)\eta^{(q+1)}_r({\rm ad}\partial_r)^{p^{q+1}}.$ Then
$$\phi^{(q+1)}(\mathrm{D_{H}}(x^{(p^{q+1}+1)\varepsilon_i)}))=0 \quad
 \mbox{for  all} \  \ i \in Y_0.$$

Suppose  $p^{q+1}+1<a_i< p^{q+2}.$ Then $a_i $ and $a_i-1 $
are not any $p$-power. Thus, by Lemmas \ref{L3.1}
and \ref{L3.2},
using induction on $a_i$, one may show that
$$\phi^{(q+1)}(\mathrm{D_{H}}(x^{a_i\varepsilon_i)}))=0 \quad \mbox{for all}\,\, i \in Y_0\, \mbox{and}\,
\ p^{q+1}+1<a_i<p^{q+2}.$$ Summarizing,
$\phi^{(q+1)}(\mathrm{D_{H}}(x^{a_i\varepsilon_i)}))=0 $ for all
$i \in Y_0$
and $1\leq a_i<p^{q+2}.$ %%%%%%%%%%%%%%%7%%%%%%%%%%%%%%%%%%%%%
Thus we construct inductively a series of derivations
$$\phi=\phi^{(0)},\phi^{(1)},\ldots,\phi^{(k)},\ldots$$
Let $k=\max\{t_1,\ldots,t_{2m}\}.$ Then
$\phi^{(k)}(\mathrm{D_{H}}(x^{a_i\varepsilon_i)}))=0 $ for all
$i\in Y_0$ and $a_i\leq \pi_i.$ Moreover, by the process of
construction , $\phi^{(k)}$ coincides  with
$$\phi-\sum_{r\in
Y_0}\sum_{s_r=1}^{t_r-1}(\lambda^{(s_r)}_r\Phi^{(s_r)}_r
+\mu^{(s_r)}_r\Theta_{r}^{s_{r}}+\eta^{(s_r)}_r({\rm
ad}\partial_r)^{p^{s_r}})
$$
on the elements $\mathrm{D_{H}}(x^{(a_i\varepsilon_i)})$ for all
$i\in Y_0$ and $a_i\leq \pi_i.$ The proof is complete.
\end{proof}
For convenience, put $\Gamma':=\sum_{r\in Y_1}x_{r}\partial_{r}.$ For the elements in $\mathcal{N},$ we have:

\begin{lemma}\label{L3.3}
Let $\phi\in {\rm Der}(\frak{N},\mathcal{W})$ be homogeneous such
that $\phi(\frak{N}_{[-1]}\oplus \frak{N}_{[0]})=0.$ Then there is
$\lambda \in \mathbb{F}$ such that $(\phi-\lambda {\rm ad}
\Gamma')(\mathcal{N})=0.$
\end{lemma}
 \begin{proof} Recall that
$\mathcal{N}=\{\mathrm{D_{H}}(x_ix^u)\ |\ i\in Y_0,\ u\in
\mathbb{B}_2\}.$ Since $\phi (\frak{N}_{[0]})=0,$ we may assume
that
\begin{equation}
\phi(\mathrm{D_{H}}(x_ix^u))=\sum_{r\in Y}f_{i,u,r}\partial_r\quad {\rm
where}\,\; f_{i,u,r} \in \Lambda(n),\;\, i\in Y_0,\;u\in
\mathbb{B}_2.\label{e6.6.31}
\end{equation}%3.7
If ${\rm zd}(\phi)$ is odd then ${\rm
zd}(\phi(\mathrm{D_{H}}(x_ix^u)))$ is even and therefore,
$f_{i,u,r}=0$ for all $r\in Y_0.$ Thus we obtain from
(\ref{e6.6.31}) that
\begin{equation}
\phi(\mathrm{D_{H}}(x_ix^u))=\sum_{r\in
Y_1}f_{i,u,r}\partial_r.\label{e6.6.32}
\end{equation}%3.8
Note that
$[\mathrm{D_{H}}(x_{i'}x_i),\mathrm{D_{H}}(x_ix^u)]=\tau(i')\mathrm{D_{H}}(x_ix^u).$
Applying $\phi$ to this equation, we obtain from (\ref{e6.6.32})
that $\phi(\mathrm{D_{H}}(x_ix^u))=0.$

Now suppose   ${\rm zd}(\phi)$ is even. Then $f_{i,u,r}=0$ for
all $r\in Y_1$ and therefore,
\begin{equation}
\phi(\mathrm{D_{H}}(x_ix^u))=\sum_{r\in
Y_0}f_{i,u,r}\partial_r.\label{e6.6.33}
\end{equation}%3.9
Obviously,
\begin{equation}
[\mathrm{D_{H}}(x^{(2\varepsilon_i)}),\phi(\mathrm{D_{H}}(x_{i'}x^u))]
=\tau(i)\phi(\mathrm{D_{H}}(x_{i}x^u)).\label{e6.6.34}
\end{equation}%3.10
From(\ref{e6.6.33}) and (\ref{e6.6.34}) we have
\begin{equation}
\phi(\mathrm{D_{H}}(x_ix^u))=f_{i,u,i'}\partial_{i'}.\label{e6.6.35}
\end{equation}%3.11
By our hypothesis that $n\geq 4,$ we can find $k,l\in Y_1\setminus
u$ with $k\not=l.$ Then
$$[\mathrm{D_{H}}(x_kx_l),\mathrm{D_{H}}(x_ix^u)] =0.$$ Applying
$\phi$ to this equation and using (\ref{e6.6.35}) we obtain that
$(x_l\partial_k-x_k\partial_l)(f_{i,u,i'})=0$ and therefore,
$x_k\partial_l(f_{i,u,i'})=0.$ This forces $f_{i,u,i'}\in
\mathbb{F}\cup\mathbb{F}\, x^{\omega}\cup\mathbb{F}\, x^u. $
Let $u=<r,s>.$ Take $k\in Y_1\setminus \{r,s\}.$ Then
$[\mathrm{D_{H}}(x_kx_r),\mathrm{D_{H}}(x_ix_kx_s)]=\mathrm{D_{H}}(x_ix^u)$
and therefore,
\begin{equation}
[\mathrm{D_{H}}(x_kx_r),\phi(\mathrm{D_{H}}(x_ix_kx_s))]
=\phi(\mathrm{D_{H}}(x_ix^u)).\label{e6.6.36}
\end{equation}%3.12
If $f_{i,u,i'}\in \mathbb{F}\cup\mathbb{F}\, x^{\omega},$ then
it is easily seen from (\ref{e6.6.36}) that
$\phi(\mathrm{D_{H}}(x_ix^u))=0.$ If $f_{i,u,i'}\in\mathbb{F}\,
x^u,$ let $f_{i,u,i'}=\lambda_{i,u,i'}x^u.$ Denote
$\lambda_{i,u}:=\lambda_{i,u,i'}$ in the following. Then we obtain
from (\ref{e6.6.36}) that $\lambda_{i,<k,s>} =\lambda_{i,<r,s>}.$
This proves that $\lambda_{i,u}$ is independent of the choice of
$u$. We denote $\lambda_{i}:=\lambda_{i,u}$ for $i\in Y_0$. Then
$\phi(\mathrm{D_{H}}(x_ix^u))=\lambda_{i}x^u\partial_{i'}$ for all $u\in
\mathbb{B}_2.$ On the other hand, for $i,j\in Y_0,$ we have
$$[\mathrm{D_{H}}(x_ix_{j'}), \mathrm{D_{H}}(x_jx^u))]=\tau(j')
(1+\delta_{i,j'})\mathrm{D_{H}}(x_ix^u).$$ Applying $\phi$ to the
equation above, one gets
$$\tau(i)\lambda_{i}=\tau(j)\lambda_{j} \quad i,j\in Y_0.$$
Let $\lambda:=\frac{1}{2}\tau(i)\lambda_{i}.$ Then
$$(\phi-\lambda{\rm ad}\Gamma')(\mathrm{D_{H}}(x_ix^u))=\lambda_{i}x^u\partial_{i'}-2\tau(i)
\lambda x^u\partial_{i'}=0.$$ The proof is complete.
\end{proof}

Now we can give a explicit description for the derivations from
$\frak{N}$ into $\mathcal{W}$ vanishing on the top of $\frak{N}.$

\begin{theorem}\label{t3.6}
Let $\phi\in {\rm Der}(\frak{N},\mathcal{W})$ be homogeneous and
 $\phi(\frak{N}_{[-1]}\oplus \frak{N}_{[0]})=0.$ Then there are
$\lambda,$ $\lambda^{(s_r)}_r,$ $ \mu^{(s_r)}_r,$ $
\eta^{(s_r)}_r\in \mathbb{F},$ where $r\in Y_0$ and $1\leq s_r <
t_r, $ such that $\mu^{(s_r)}_r\eta^{(s_r)}_r=0$ and
$$\phi=\lambda{\rm ad}\Gamma'+\sum_{r=1}^{2m}
\sum_{s_r=1}^{t_r-1}(\lambda^{(s_r)}_r\Phi^{(s_r)}_r
+\mu^{(s_r)}_r \Theta_{r}^{s_{r}}+\eta^{(s_r)}_r({\rm
ad}\partial_r)^{p^{s_r}}).$$
In particular, if $n$ is odd then $\phi=\lambda{\rm ad}\Gamma'+\sum_{r=1}^{2m}
\sum_{s_r=1}^{t_r-1}\eta^{(s_r)}_r({\rm
ad}\partial_r)^{p^{s_r}}.$
\end{theorem}
 \begin{proof} By Lemmas \ref{L3.4}
 and \ref{L3.5},
there are $\lambda^{(s_r)}_r,$ $ \mu^{(s_r)}_r,$ $
\eta^{(s_r)}_r\in \mathbb{F}$ such that
$$\widetilde{\phi}:=\phi-\sum_{r=1}^{2m}
\sum_{s_r=1}^{t_r-1}(\lambda^{(s_r)}_r\Phi^{(s_r)}_r
+\mu^{(s_r)}_r\Theta_{r}^{s_{r}}+\eta^{(s_r)}_r({\rm
{ad}}\partial_r)^{p^{s_r}})
$$ vanishes on $\mathcal{M}$. By Lemma \ref{L3.3}
there is $\lambda\in \mathbb{F}$ such that
$\varphi:=\widetilde{\phi}-\lambda{\rm ad}\Gamma'$ vanishes on
$\mathcal{N}.$ It is easy to see that $\varphi$ vanishes on
$\mathcal{M}$ and Theorem \ref{t1.7}
ensures that $\varphi=0.$ The remaining is clear.
 \end{proof}

It is our present aim to reduce certain derivations to be vanishing the top of $\frak{N}.$
 \begin{lemma} \label{P2.1.6} (see \cite[Proposition 2.1.6]{lz1})
Let $\mathcal{L}$
  be a $\mathbb{Z}$-graded subalgebra of $\mathcal{W}$ such that
  $\mathcal{L}_{-1}=\mathcal{W}_{-1}.$
  If $\phi \in \mathrm{Der}_{t}(\mathcal{L},\mathcal{W})$
  where $t:=\mathrm{zd}(\phi)\geq 0,$
  then there exists $E\in \mathcal{W}_{t}$
  such that
  $$\big(\phi -\mathrm{ad}E\big)(\mathcal{L}_{-1})=0.$$
  \end{lemma}

  In view of Lemma  \ref{P2.1.6} and Theorem \ref{t3.6},
 it suffices to reduce the homogeneous derivations vanishing on
$\frak{N}_{[-1]}$ to be vanishing on the top of $\frak{N}.$

\begin{lemma}\label{a6.6.14} Let $\phi\in {\rm Der}(\frak{N},\mathcal{W})$ be
$\mathbb{Z}$-homogeneous with odd degree such that
$\phi(\frak{N}_{[-1]})=0.$ Then there is $D\in\mathcal{W}$ such
that $(\phi-{\rm ad}D)(\frak{N}_{[-1]}\oplus\frak{N}_{[0]} )=0.$
\end{lemma}
 \begin{proof} Since $\phi(\frak{N}_{[-1]})=0$ and ${\rm
zd}(\phi)$ is odd, by Lemma \ref{bc2.2}, $\phi(\frak{N}_{[0]})\in O(\mathcal{G}).$ For
$1\leq i\leq m,$ one may assume that
\begin{equation}
\phi(\mathrm{D_{H}}(x_ix_{i'}))=\sum _{r\in
Y_0}f_{i,r}\partial_{r}\quad{\rm where}\  f_{i,r}\in
\Lambda(n).\label{e6.6.37}
\end{equation}%4.1
Let $ 1\leq i\not= j\leq m.$ Then $[\mathrm{D_{H}}(x_ix_{i'}),\
\mathrm{D_{H}}(x_jx_{j'})]=0 $ and therefore,
\begin{equation} [\phi(\mathrm{D_{H}}(x_ix_{i'})),\ \mathrm{D_{H}}(x_jx_{j'})]
=[\phi(\mathrm{D_{H}}(x_jx_{j'})),\
\mathrm{D_{H}}(x_ix_{i'})].\label{e6.6.38}
\end{equation}%4.2
One may easily deduce from (\ref{e6.6.37}) and (\ref{e6.6.38})
that $f_{ij}=f_{ij'}=0.$ Thus,
$$\phi(\mathrm{D_{H}}(x_ix_{i'}))=f_{ii}\partial_i+f_{ii'}\partial_{i'}\quad {\rm for}\
 1\leq i\leq m.$$
Put $\psi:=\phi - {\rm
ad}(\sum^{m}_{r=1}(f_{rr'}\partial_{r'}-f_{rr}\partial_{r} )). $ Then
$\psi(\mathrm{D_{H}}(x_ix_{i'}))=0$ for $i\in Y_0.$ We propose to
show that
$$\psi(\mathrm{D_{H}}(x_ix_{j}))=0\quad\mbox{for all}\ i,j\in Y_0.
$$
First, we consider the case $j\not= i,i'.$ Clearly,
$[\mathrm{D_{H}}(x_ix_{i'}),\
\mathrm{D_{H}}(x_ix_{j})]=\tau(i')\mathrm{D_{H}}(x_ix_{j}).$
Applying $\psi$ to this equation, we have
$$[\mathrm{D_{H}}(x_ix_{i'}),\ \psi(\mathrm{D_{H}}(x_ix_{j}))]
=\tau(i')\psi(\mathrm{D_{H}}(x_ix_{j})).
$$
Similarly,
$$[\mathrm{D_{H}}(x_jx_{j'}),\ \psi(\mathrm{D_{H}}(x_ix_{j}))]
=\tau(j')\psi(\mathrm{D_{H}}(x_ix_{j})).$$ Noticing
$\psi(\mathrm{D_{H}}(x_ix_{j}))\in O(\mathcal{G})$ and $ j\not= i,i',$ one gets from the above two equations that
$\psi(\mathrm{D_{H}}(x_ix_{j})=0.$

Next, we consider $\psi(\mathrm{D_{H}}(x^{(2\varepsilon_i)})).$
Clearly,
$$[\mathrm{D_{H}}(x_ix_{i'}), \ \mathrm{D_{H}}(x^{(2\varepsilon_i)})]
=2\tau(i')\mathrm{D_{H}}(x^{(2\varepsilon_i)}).
$$
Assume that $\psi(\mathrm{D_{H}}(x^{(2\varepsilon_i)}))=\sum_{r\in
Y_0}g_{ir}\partial_r,\ g_{ir}\in \mathbb{F}.$ Applying $\psi$ to the
equation above, we have $g_{ir}=0$ for $r\not= i,i'$ and
$3\tau(i')g_{ii}=0,\tau(i)g_{ii'}=0.$ Thus
$\psi(\mathrm{D_{H}}(x^{(2\varepsilon_i)}))=0.$

For $k,l\in Y_1,$ applying $\psi$ to the equation
$[\mathrm{D_{H}}(x_ix_{i'}),\ \mathrm{D_{H}}(x_kx_{l})]=0$
 for $i\in Y_0,$
 one may easily show that
 $\psi(\mathrm{D_{H}}(x_kx_{l}))=0.$
 Summarizing, we have shown that
 $\psi(\frak{N}_{[0]})=0$ and then
  $\psi(\frak{N}_{[-1]}\oplus\frak{N}_{[0]})=0.$
  \end{proof}

  Put ${\rm
Der}^-(\frak{N},\mathcal{W}):=\oplus_{r\leq -1}{\rm
Der}_r(\frak{N},\mathcal{W}). $
\begin{theorem}\label{t3.8} If $n$ is even, then
$${\rm Der}^-(\frak{N},\mathcal{W})={\rm ad}\mathcal{W}_{[-1]}
+\sum_{1\leq r\leq 2m}\sum_{\stackrel{1\leq s_r\leq
t_r-1}{n-p^{s_r}<0}}
(\mathbb{F}\,\Phi^{(s_r)}_r+\mathbb{F}\, \Theta_{r}^{s_{r}})
+\sum_{r\in Y_0}\sum_{1\leq s_r\leq t_r-1} \mathbb{F}\,({\rm
ad}\partial_r)^{p^{s_r}}.$$
If $n$ is odd, then ${\rm Der}^-(\frak{N},\mathcal{W})={\rm ad}\mathcal{W}_{[-1]}
+\sum_{r\in Y_0}\sum_{1\leq s_r\leq t_r-1} \mathbb{F}\,({\rm
ad}\partial_r)^{p^{s_r}}.$
\end{theorem}
 \begin{proof} Suppose   $n$ is even. Let $\phi\in {\rm
Der}(\frak{N},\mathcal{W})$ be homogeneous. If ${\rm zd}(\phi)\leq
-2,$ then Theorem \ref{t3.6}
 shows that
$$
\phi\in \sum_{1\leq r\leq 2m }\sum_{\stackrel{1\leq s_r\leq
t_r-1}{n-p^{s_r}<0}}\Big(\mathbb{F}\,\Phi^{(s_r)}_r
+\mathbb{F}\,\Theta_{r}^{s_{r}} +\sum_{r\in Y_0}\sum_{1\leq
s_r\leq t_r-1} \mathbb{F}\,({\rm ad}\partial_r)^{p^{s_r}}\Big).
$$ Now assume that ${\rm
zd}(\phi)=-1.$  By Lemma \ref{a6.6.14},
 there is
$D\in\mathcal{W}$ such that $(\phi-{\rm
ad}D)(\frak{N}_{[-1]}\oplus\frak{N}_{[0]} )=0.$ Again Theorem
\ref{t3.6}
 entails that
$$ \phi\in{\rm ad}\mathcal{W}+\mathbb{F}\, \mathrm{ad}\Gamma'
+\sum_{1\leq r\leq 2m}\sum_{ 1\leq s_r\leq t_r-1}
(\mathbb{F}\,\Phi^{(s_r)}_r+\mathbb{F}\, \Theta_{r}^{s_{r}})
+\sum_{r\in Y_0}\sum_{1\leq s_r\leq t_r-1} \mathbb{F}\,({\rm
ad}\partial_r)^{p^{s_r}}.
$$
Since $\phi$ is $\mathbb{Z}$-degree $-1,$ one may easily obtain
the desired result.  In the case that $n$ is odd, the argument is similar.
\end{proof}

As an direct application of Theorem \ref{t3.6}
 and Lemma \ref{a6.6.14},
  we can  determine all the
  derivations  of
odd $\mathbb{Z}$-degree:

\begin{theorem}\label{t3.9}   If $n$ is even, then
$$\sum_{k\, \mathrm{odd}} \mathrm{Der}_{[k]}(\frak{N},\mathcal{W})=
\sum_{i\, \mathrm{odd}} \mathrm{ad}\mathcal{W}_{[i]}+\sum_{1\leq r\leq 2m}\sum_{1\leq
s_r\leq t_r-1}
(\mathbb{F}\,\Phi^{(s_r)}_r+\mathbb{F}\, \Theta_{r}^{s_{r}})
+\sum_{r\in Y_0}\sum_{1\leq s_r\leq t_r-1} \mathbb{F}\,({\rm
ad}\partial_r)^{p^{s_r}}.$$
If $n$ is odd, then
$$\sum_{k\, \mathrm{odd}} \mathrm{Der}_{[k]}(\frak{N},\mathcal{W})=
\sum_{i\, \mathrm{odd}} \mathrm{ad}\mathcal{W}_{[i]}
+\sum_{r\in Y_0}\sum_{1\leq s_r\leq t_r-1} \mathbb{F}\,({\rm
ad}\partial_r)^{p^{s_r}}.$$
\end{theorem}

Let us consider the derivations  of even $\mathbb{Z}$-degree.

\begin{proposition}\label{p3.10}
Let $\phi\in {\rm Der}(\frak{N},\mathcal{W})$ be  of
  even $\mathbb{Z}$-degree such that
$\phi(\frak{N}_{[-1]})=0$. Then there are
$\lambda_{r}\in\mathbb{F}\ (1\leq r\leq m)$ such that
$$
\Big(\phi-\sum_{1\leq r\leq m}\lambda_{r}\Psi^{(r)}\Big)
(\mathrm{D_{H}}(x_ix_{j}))=0\quad \mbox{for all}\  i,j\in Y_0.
$$
\end{proposition}
 \begin{proof}  Since ${\rm zd}(\phi)$ is even and
$\phi(\frak{N}_{[-1]})=0,$   Lemma \ref{bc2.2} ensures that $ \phi(\frak{N}_{[0]})\subset
E(\mathcal{G}).$ Thus one may assume that for $1\leq i\leq m,$
$$\phi(\mathrm{D_{H}}(x_ix_{i'}))
=\sum_{r\in Y_1}f_{ir}\partial_r\quad \mbox{where}\ f_{ir}\in
\Lambda(n).$$ Applying $\phi$ to
$[\mathrm{D_{H}}(x_kx_{l}), \mathrm{D_{H}}(x_ix_{i'})]=0$ for
$k,l\in Y_1,$ we have
$$[\mathrm{D_{H}}(x_kx_{l}),\
\phi(\mathrm{D_{H}}(x_ix_{i'}))]=0
$$
since $\phi(\mathrm{D_{H}}(x_kx_{l}))\in E(\mathcal{G}).$ Thus,
 $[x_l\partial_k-x_k\partial_l,\sum_{i\in Y_1}f_{ir}\partial_r]=0$  and therefore,
\begin{equation}
\sum _{r\in
Y_1}(x_l\partial_k-x_k\partial_l)(f_{ir})\partial_{r}-f_{il}\partial_{k}+f_{ik}\partial_{l}=0.\label{e6.6.39}%4.3
\end{equation}
This implies that
\begin{equation}
(x_l\partial_k-x_k\partial_l)(f_{ir})=0\quad {\rm whenever}\  k,l\not=r;
\label{e6.6.40}
\end{equation}%4.4
\begin{equation}
(x_l\partial_k-x_k\partial_l)(f_{ik})-f_{il}=0.\label{e6.6.41}
\end{equation}%4.5
By (\ref{e6.6.40}), one may deduce that $f_{ir}\in \mathbb{F}\,
x_r\cup\mathbb{F}\, x^{\omega}\cup\mathbb{F}\,
x^{\omega-<r>}.$ Assume that $|\omega|=n$ is odd and
$0\not=f_{ir}\in \mathbb{F}\, x^{\omega}.$ Then one may reach a
contradiction from (\ref{e6.6.39}). It remains to consider the
following two cases, since $\phi$ is homogeneous.\newline

 \noindent\textit{Case} (i): $f_{ir}\in \mathbb{F}x_r$ for every
$r\in Y_1.$ Assume that
$$\phi(\mathrm{D_{H}}(x_ix_{i'}))=\sum _{r\in
Y_1}\mu_{ir}x_r\partial_r,\quad \mu_{ir}\in \mathbb{F}.$$ From
(\ref{e6.6.41}) we have $\mu_{ir}=\mu_{is}$ for all $r,s\in Y_1$.
Denote $\mu_{i}:=\mu_{ir}$ for $r\in Y_1$. Thus,
$$\phi(\mathrm{D_{H}}(x_ix_{i'}))=\mu_{i}\Gamma',\quad r\in Y_0.$$
Recall our assumption that $n\geq4.$ Let $u\in \mathbb{B}_4.$
Clearly, $\phi(\mathrm{D_{H}}(x^u))\in E(\mathcal{G}).$ Applying
$\phi$ to the equation $[\mathrm{D_{H}}(x_ix_{i'}),\
\mathrm{D_{H}}(x^{u})]=0,$ we have
$$
2\mu_i\mathrm{D_{H}}(x^{u})=[\mu_i\Gamma',\
\mathrm{D_{H}}(x^{u})]= [\phi(\mathrm{D_{H}}(x_ix_{i'})),\
\mathrm{D_{H}}(x^{u})]=0.
$$
Consequently, $\mu_i=0.$ Hence $\phi(\mathrm{D_{H}}(x_ix_{i'}))=0$
for $i\in Y_0.$ For $i,j\in Y_0$ with $j\not=i',$ we have
$[\mathrm{D_{H}}(x_ix_{i'}),\
\mathrm{D_{H}}(x_ix_{j})]=\tau(i')(1+\delta_{ij})\mathrm{D_{H}}(x_ix_{j}).
 $ It follows that $\phi(\mathrm{D_{H}}(x_ix_{j}))=0 $ for $i,j\in
 Y_0$ with $j\not= i'.$ Hence
 $$\phi(\mathrm{D_{H}}(x_ix_{j}))=0\quad {\rm for}\ i,j\in Y_0.$$

 \noindent\textit{Case} (ii): $f_{ir}\in \mathbb{F}\,
x^{\omega-<r>}$ for all $r\in Y_1.$ Assume that
$$\phi(\mathrm{D_{H}}(x_ix_{i'}))=\sum _{r\in
Y_1}\lambda_{ir}x^{\omega-<r>}\partial_r,\quad \lambda_{ir}\in
\mathbb{F}.$$ For $k,l\in Y_1$ with $k\not=l,$ we have
$[\phi(\mathrm{D_{H}}(x_ix_{i'})),\ \mathrm{D_{H}}(x_kx_{l})]=0.$
A direct  computation shows that
$(-1)^k\lambda_{ik}=(-1)^l\lambda_{il}.$ Denote
$\lambda_{i}:=(-1)^k\lambda_{ik}$ for $k\in Y_1.$ Then
$$\phi(\mathrm{D_{H}}(x_ix_{i'}))=\lambda_{i}\mathrm{D_{H}}(x^{\omega}).$$
Put $\psi:=\phi-\sum_{1\leq r\leq m} \lambda_{r}\Psi^{(r)}.$ Then
$\psi(\mathrm{D_{H}}(x_ix_{i'}))=0.$ Just as in Case (i), we have
$\psi(\mathrm{D_{H}}(x_ix_{j}))=0 $ for all $i,j\in Y_0.$ The
proof is complete.
\end{proof}

\begin{remark}\label{r3.11} In view of Lemma \ref{P2.1.6},
Proposition \ref{p3.10},
 Theorems \ref{t3.6},
 \ref{t3.8}  and
  and \ref{t3.9},
  in order to determine the
derivation space of $\frak{N}$ to $\mathcal{W},$ it  suffices to
reduce every even $\mathbb{Z}$-degree  derivation  in
$\mathrm{Der}(\frak{N},\mathcal{W})$   vanishing on
$\frak{N}_{[-1]}$ to be the one vanishing on
$\big\{\mathrm{D_H}(x_{k}x_{l})\mid k, l\in Y_1\big\} $ modulo a
suitable known derivation (such as a linear combination of an
inner derivation and some exceptional derivations). However, in
contrast to the situation of $HO$ or $S$ (see \cite{lz1, lhs}), it
seems that we must find a new approach other than the method of
canonical torus used there.  We believe that this problem will be settled in
future.
\end{remark}
Using Theorems \ref{t1.4}, \ref{t3.6}, \ref{t3.8}, and \ref{t3.9}, we obtain the corresponding
results on the even part $\mathcal{H}$ of the finite-dimensional Hamiltonian superalgebra:
\begin{corollary}
Let $\phi\in {\rm Der}(\mathcal{H},\mathcal{W})$ and
 $\phi(\mathcal{H}_{[-1]}\oplus \mathcal{H}_{[0]})=0.$ Then there are
$\lambda,$ $\mu,$ $\lambda^{(s_r)}_r,$ $ \mu^{(s_r)}_r,$ $
\eta^{(s_r)}_r\in \mathbb{F},$ where $r\in Y_0$ and $1\leq s_r <
t_r, $ such that $\mu^{(s_r)}_r\eta^{(s_r)}_r=0$ and
$$\phi=\Gamma_{\lambda} +\mu({\rm ad}\Gamma')+\sum_{r=1}^{2m}
\sum_{s_r=1}^{t_r-1}(\lambda^{(s_r)}_r\Phi^{(s_r)}_r
+\mu^{(s_r)}_r \Theta_{r}^{s_{r}}+\eta^{(s_r)}_r({\rm
ad}\partial_r)^{p^{s_r}}).$$
\end{corollary}

\begin{corollary} If $n$ is even, then
$${\rm Der}^-(\mathcal{H},\mathcal{W})={\rm ad}\mathcal{W}_{[-1]}
+\mathbb{F}\, (\delta_{|\pi|>|\omega|}\Gamma_{1}) +\sum_{1\leq r\leq
2m}\sum_{\stackrel{1\leq s_r\leq t_r-1}{n-p^{s_r}<0}}
(\mathbb{F}\,\Phi^{(s_r)}_r+\mathbb{F}\, \Theta_{r}^{s_{r}})
+\sum_{r\in Y_0}\sum_{1\leq s_r\leq t_r-1} \mathbb{F}\,({\rm
ad}\partial_r)^{p^{s_r}},$$
where $\Gamma_{1}$ is defined as in Remark \ref{r1.6}.
If $n$ is odd, then
$${\rm Der}^-(\mathcal{H},\mathcal{W})={\rm ad}\mathcal{W}_{[-1]}
+\sum_{r\in Y_0}\sum_{1\leq s_r\leq t_r-1} \mathbb{F}\,({\rm
ad}\partial_r)^{p^{s_r}}.$$
\end{corollary}
\begin{corollary}
   If $n$ is even, then
$$ \sum_{k\, \mathrm{odd}}\mathrm{Der}_{[k]}(\mathcal{H},\mathcal{W})=
\sum_{i\, \mathrm{odd}}\mathrm{ad}\mathcal{W}_{[i]}+\sum_{1\leq r\leq
2m}\sum_{1\leq s_r\leq t_r-1}
(\mathbb{F}\,\Phi^{(s_r)}_r+\mathbb{F}\, \Theta_{r}^{s_{r}})
+\sum_{r\in Y_0}\sum_{1\leq s_r\leq t_r-1} \mathbb{F}\,({\rm
ad}\partial_r)^{p^{s_r}}.$$
If $n$ is odd, then
$$ \sum_{k\, \mathrm{odd}}\mathrm{Der}_{[k]}(\mathcal{H},\mathcal{W})=
\sum_{i\, \mathrm{odd}}\mathrm{ad}\mathcal{W}_{[i]}+\sum_{r\in Y_0}\sum_{1\leq s_r\leq t_r-1} \mathbb{F}\,({\rm
ad}\partial_r)^{p^{s_r}}.$$
\end{corollary}


\begin{thebibliography}{99}

\bibitem{Ce} M. J. Celousov. Derivations of Lie algebras of
Cartan type. \textit{Izv. Vyss}. \textit{Ucebn}. \textit{Zaved.
Mathematika}  \textbf{98} (1970): 126--134 (in Russian).

\bibitem{Kac2} V. G. Kac. Lie superalgebras. \textit{Adv.
Math.}  \textbf{26} (1977): 8--96.

\bibitem{Kac4} V. G. Kac. Classification of supersymmetries.
                     \textit{Proc. ICM}  Vol. \textbf{1} (2002): 319--344.

\bibitem{KL} Yu. Kochetkov. and  D. Leites. Simple Lie algebras in
               characteristic 2 recovered from superalgebras and on the notion of
             a simple finite group. \textit{Contemp. Math.} Vol. \textbf{131},
          Part 2, \textit{Amer. Math. Soc., Providence, RI}, 1992.
\bibitem{lg}  W.-D. Liu and B.-L. Guan. Derivations  from  the even parts into the odd
          parts for Lie superalgebras $W$ and $S$. \textit{arXive}: math. R.A/0509177.

\bibitem{lhs} W.-D. Liu, X.-Y. Hua, and Y.-C. Su. Derivations for the even part of the odd Hamiltonian superalgebra
in positive characteristic.

\bibitem{lz1}  W.-D. Liu and  Y.-Z. Zhang.
 Derivations for the even parts of modular Lie superalgebras $W$
 and $S$ of Cartan type. \textit{arXive}: math. R.A/0507521.

\bibitem{LZ1} W.-D. Liu,  Y.-Z. Zhang and  X.-L. Wang. The
derivation algebra of the Cartan-type Lie superalgebra $HO$,
\textit{J. Algebra} \textbf{273}(1) (2004), 176--205.

\bibitem{MZ} F.-M. Ma and  Q.-C. Zhang. Derivation algebra of modular
Lie superalgebra
    $K$ of Cartan type. \textit{J. Math. (PRC)}.
    \textbf{20}(4) (2000): 431--435 (in Chinese).

\bibitem{P} V. M. Petrogradski. Identities in the enveloping algebras for
 modular Lie superalgebras. \textit{J. Algebra}  \textbf{145} (1992): 1--21.

\bibitem{Sch1}M. Scheunert. ``\textit{Theory of Lie
superalgebras}", in: Lecture Notes in Math.
 \textbf{716}, Springer-Verlay, 1979.

\bibitem{st} H. Strade. \textit{Simple Lie algebras over fields of positive characteristic,
I. Structure theory.} Walter de Gruyter, Berlin and New York, 2004.


 \bibitem{SF}
H. Strade and R. Farnsteiner. ``\textit{Modular Lie algebras and
their representations}'', in:  Monogr. Textbooks Pure Appl. Math.
Vol. 116, Dekker, 1988.



\bibitem{WZ3} Y. Wang and Y.-Z. Zhang. Derivation algebra
$\mathrm{Der}(H)$  and central extensions of Lie superalgbras.
\textit{Commun. Algebra}  \textbf{32} (2004): 4117--4131.

 \bibitem{ZF}  Y.-Z. Zhang and H.-C. Fu.  Finite-dimensional
 Hamiltonian Lie superalgebras. \textit{Commun. Algebra}
 \textbf{30}(6) (2002): 2651--2673.

\bibitem{Zh1} Y.-Z. Zhang. Finite-dimensonal Lie superalgebras of
               Cartan-type over fields of prime characteristic.
              \textit{Chin. Sci.  Bull.}  \textbf{42} (1997): 720--724.

\bibitem{ZZ1} Q.-C. Zhang and Y.-Z. Zhang.  Derivation algebras
of the modular
    Lie superalgebras $W$ and $S$ of Cartan type. \textit{Acta Math. Sci.}
    \textbf{20}(1) (2000): 137--144.
 \end{thebibliography}
\end{document}